\documentclass[letterpaper, paper,11pt]{AAS}	

\usepackage{bm}
\usepackage{amsmath}
\usepackage{subcaption}
\usepackage[colorlinks=true, pdfstartview=FitV, linkcolor=black, citecolor= black, urlcolor= black]{hyperref}
\usepackage{overcite}
\usepackage{footnpag}			      	
\usepackage{graphicx}
\usepackage{algpseudocode}
\usepackage{algorithm}
\usepackage{bm}
\usepackage{graphicx}
\usepackage{amsfonts}
\usepackage{algorithm}
\usepackage{algpseudocode}

\PaperNumber{24-322}

\begin{document}

\title{CO-OPTIMIZATION OF SPACECRAFT AND LOW-THRUST TRAJECTORY WITH DIRECT METHODS}

\author{Keziban Saloglu\thanks{PhD Student, Department of Aerospace Engineering, Auburn University, Auburn, AL 36849, USA, AAS Member.} 
\ and Ehsan Taheri\thanks{Assisstant Professor, Department of Aerospace Engineering, Auburn University, Auburn, AL 36849, USA, AAS Member.}
}

\maketitle{}

\begin{abstract}
Solar-powered electric propulsion systems can operate in multiple modes and their operation is coupled to the power generated by solar arrays. However, the power produced by the solar arrays is a function of the solar array size and heliocentric distance to the Sun, which also depends on the to-be-optimized trajectory. The optimization of spacecraft solar array size, thruster modes, and trajectory can be performed simultaneously, capitalizing on the inherent couplings. In this work, we co-optimize the spacecraft's solar array size, thruster modes, and trajectory using a direct optimization, which allows for maximizing the net delivered mass. 
A particular challenge arises due to the existence of discrete operating modes. We proposed a method for smoothly selecting optimal operation modes among a set of finite possible modes.
The utility of the proposed method is demonstrated successfully by solving a benchmark problem, i.e., the Earth to Comet 67P fixed-time rendezvous problem. The results indicate that utilizing multiple modes increases the net useful mass compared to a single mode and leads to a smaller solar array size for the spacecraft. 
\end{abstract}

\section{Introduction}
Solar electric propulsion (SEP) systems offer high efficiency with respect to propellant consumption due to their high specific impulse values. Many space missions have  successfully utilized SEP thrusters with different scientific mission goals\cite{kawaguchi2008hayabusa, rayman2000results, russell2012dawn, oh2019development}. In low-thrust trajectory optimization with SEP systems, and for a more realistic mission design, power constraints have to be taken into consideration. The power subsystem of the spacecraft highly influences the propulsion subsystem and the profile of the trajectory, because the SEP thrusters require solar power to operate and based on the available power and optimality can switch between different operation modes.

Typically, electric thrusters can be configured to operate at multiple discrete operation modes, which are characterized by specific input power, thrust level, and mass flow rate values~\cite{patterson2007next}. These operation modes can be configured during the spacecraft design phase. Since the discrete operation modes of the thrusters change with different power inputs, the sizing of the solar arrays affects the performance of the thrusters \cite{quarta_minimum-time_2011,petukhov_joint_2019}. One common approach to model thrusters is to consider the thrust level and the mass flow rate as quartic (or quintic) polynomials, with input power as the independent variable ~\cite{sauer_jr_modeling_1978,rayman_design_2002,TAHERI2020166}. For example, multimode thruster models are used for minimum-time low-thrust trajectory design in Ref.~\citenum{quarta_minimum-time_2011} and maximization of net useful payload mass in Ref.~\citenum{arya_electric_2021}.

The power produced by the solar array is a function of both its size and the heliocentric distance of the spacecraft. Due to the inherent couplings between the trajectory, power, and propulsion subsystems, the operational modes of the thruster, trajectory, and solar array size can be optimized during the mission design process by formulating and solving a co-optimization problem~\cite{arya_electric_2021}. The multimode thruster model coupled with the trajectory optimization task results in a more realistic trajectory by incorporating the power dependence and makes it possible to determine the ``optimal'' modes of operation of a thruster. Then, the configuration of the SEP thrusters with the required operation modes can be achieved during the spacecraft design process. In addition, co-optimizing the solar array size along with the trajectory offers an improved mass budget estimation~\cite{chi2018power}. However, solving the resulting optimal control problems is challenging due to 1) the discrete nature of operation modes (which turns the optimization problem into a mixed-integer programming problem), 2) discontinuities due to the presence of an unknown number of control switches, 3) the presence of discontinuous power constraints, 4) the introduction of additional design variables (e.g., the beginning-of-life power for sizing solar arrays), and 5) the inherent discipline-level nonlinear couplings. 

The low-thrust trajectory optimization problems are typically solved using direct or indirect optimization methods~\cite{betts_survey_1998,trelat2012optimal}. The indirect formalism of optimal control theory \cite{bryson2018applied} is used for low-thrust trajectory optimization, including power constraints with solar array sizing and multimode thruster models. An interesting problem is to maximize the useful mass of the spacecraft, defined as the objective of the optimization problem~\cite{petukhov2019joint,arya_composite_2021,arya_electric_2021}. When the net delivered useful mass is becomes the objective, the solar array mass appears in the  problem's cost, which is characterized by the beginning-of-life power value of the solar arrays. To size the solar arrays, one can optimize the beginning-of-life power, which can be used to determine the mass of the solar arrays. 

To overcome the numerical difficulties due to the discontinuous mixed-integer optimal control problem, the composite smooth control (CSC) method is applied to perform smoothing on the selection of the operation modes of the thruster~\cite{taheri2020novel}. The CSC method is also utilized for co-optimizing spacecraft trajectory and the propulsion system, considering variable specific impulse model for the thruster~\cite{arya_composite_2021}. A generic homotopic smoothing method is introduced to handle discrete eclipse and power constraints, where the continuation process is embedded with the variation on the control bounds~\cite{xiao_generic_2024}. 
Solar array sizing of spacecraft along with the trajectory with a multimode discrete thruster model is implemented using indirect methods~\cite{arya_electric_2021}. The joint optimization of the spacecraft trajectory is performed in Ref.~\citenum{petukhov_joint_2019} in which design variables such as the beginning-of-life power of the solar arrays, exhaust velocity of the thrusters, and the maximum heliocentric distance at which the maximum power is generated are taken into account. Minimum-time and minimum-energy trajectory optimization problems are solved by considering a power constraint with switching detection~\cite{wang_indirect_2022}. Incorporation of duty-cycle constraints (that give rise to many discontinuous events) is also considered within the indirect optimization methods \cite{nurre_duty-cycle-aware_2023}. 

On the other hand, the direct optimization methods first discretize the optimal control problems and then solve the obtained nonlinear programming (NLP) problems. The presence of multiple modes of operation poses a challenge due to the discrete nature of the decision variables. However, we have already proposed a method for solving trajectory optimization with a discrete throttle table for minimum-fuel formulation and without considering co-optimization \cite{salogluMultimode}. Evolutionary optimization algorithms are used to generate the global optimal solutions~\cite{conway_survey_2012}. One application of evolutionary algorithm based trajectory optimization is demonstrated with the Evolutionary Mission Trajectory Generator (EMTG) tool for automated global trajectory optimization of interplanetary missions ~\cite{englander2012automated,englander2017automated}. They provide three different thruster models, the first model uses a custom polynomial such that over- or under-performance is avoided~\cite{knittel_improved_2017}. The second model is the stepped model, characterized by input power. The stepped model approximates the Heaviside function discrete mode switches with a logistic function. The third model, which is a two-dimensional stepped model, allows the optimizer to select modes from the full thruster grid, characterized by the power and voltage inputs.

In this paper, we investigate the co-optimization of spacecraft trajectory, solar arrays, and the thruster modes with a multimode thruster model, which implements a higher-fidelity power model compared to our previous work~\cite{salogluMultimode}. In particular, we consider a direct optimization method due to the flexibility of direct methods in handling various types of constraints.
This paper's main contribution is co-optimizing the spacecraft solar array size, trajectory, and thruster operation modes with a direct method while considering a discrete thruster model. In our previous study, we proposed a method to address the discontinuous nature of multiple modes of the thruster while solving the problem using the direct formalism of optimal control~\cite{salogluMultimode}. This work uses a more realistic power model for the optimization problem compared to our previous study. In addition to considering the multimode thruster model, we take into account the mass of the solar arrays. We introduce two additional control variables. One is for smooth mode selection of the thrusters and the second is the beginning-of-life power of the solar arrays ($P_{\text{BL}}$). The mass of the solar arrays is approximated as proportional to the $P_{\text{BL}}$. By considering a breakdown of the spacecraft total mass, we consider the problem of maximizing the useful mass, i.e., the net delivered payload mass. 

The remainder of the paper is organized as follows. In the Optimal Control Problem Formulation section, the equations of motion of the spacecraft are given. Details of the power and thruster models are explained. The smooth mode-selection algorithm, initial guess generation, and the continuation approach are explained. In the Numerical Results section, obtained solutions for a time-fixed rendezvous problem from Earth to Comet 67P/Churyumov–Gerasimenko (Earth-67P) are presented. Finally, the concluding remarks are given in the Conclusion section.  
\section{Optimal Control Problem Formulation}
The equations of motion of the spacecraft are modeled using the set of Modified Equinoctial Elements (MEEs). It is established in the literature that the set of MEEs offer numerical advantages when solving low-thrust optimization problems~\cite{taheri2016enhanced,junkins_exploration_2019}. Let the state vector of the spacecraft $\bm{x} = [p, f, g, h, k, l]^{\top}$ and the spacecraft mass, $m$, then the equations of motion can be written as,
\begin{equation}
\dot{\boldsymbol{x}}(t)=\boldsymbol{f}(\boldsymbol{x}, \boldsymbol{a})=\boldsymbol{A}_{\text{MEE}}(\boldsymbol{x}, t)+\mathbb{B}_\text{MEE}(\boldsymbol{x}, t) \boldsymbol{a},
\label{eq:dyn}
\end{equation}
where $\boldsymbol{A}_{\text{MEE}} \in \mathbb{R}^{6\times 1}$ is the unforced part of the dynamics, $\mathbb{B}_\text{MEE} \in \mathbb{R}^{6\times 6}$ is the control-influence matrix, where the explicit relations for $\boldsymbol{A}_{\text{MEE}}$ and $\mathbb{B}_\text{MEE}$ can be found in Ref.~\citenum{junkins_exploration_2019}. The acceleration vector due to the propulsion system is denoted with $\bm{a} \in \mathbb{R}^{3 \times 1}$ and is parameterized as, 
\begin{equation}
    \bm{a} = \frac{T(n)}{m}\hat{\bm{\alpha}},
    \label{eq:acc}
\end{equation}
where $n$ is the current operation mode of the thrusters, $n \in \{1,\cdots,N_\text{mode}\}$ with $N_\text{mode}$ being the total number of discrete operation modes. Therefore, the thrust level is a function of the discrete operation mode, $n$, and $\hat{\bm{\alpha}}$ is the unit vector of thrust steering direction (i.e., $\|\hat{\bm{\alpha}}\|=1$). Similarly, the spacecraft's mass time rate of change is parameterized in terms of the operation mode as,
\begin{equation}
    \dot m = -\dot{m}(n).
    \label{eq:mdot}
\end{equation}

The resulting optimal control problem is a mixed-integer programming problem where the thrust steering is the continuous control and the current operation mode, $n$, is a discrete variable. The optimal discrete mode has to be determined at each time step. To formulate the optimal control problem as an NLP problem, an additional control input, engine power control, $P_{\text{E}}$, is introduced to facilitate the selection of the operation mode through a power-based smooth blending of the operation modes. Then, the augmented control vector can be defined as $\bm{u} = [\hat{\bm{\alpha}}^{\top}, P_{\text E}]^{\top} \in \mathbb{R}^4$. 
The current operation mode is determined by defining a function that chooses the modes smoothly given the engine power and the vector associated with the total selected mode,  $\bm{P}_{\text{sel}} \in \mathbb{R}^{N_\text{mode}}$ as, 
\begin{equation}
    n = f_{\text{mode}}(P_{\text E}, \bm{P}_{\text{sel}}).
    \label{eq:fmode}
\end{equation}

We note that $\bm{P}_{\text{sel}}$ is based on the number of modes that are selected from a throttle table. For instance, later we will consider an SPT-140 Hall thruster in formulating several optimal control problems. However, we will only select a sub-set of all modes, and those selected modes will constitute the elements of the $\bm{P}_{\text{sel}}$ vector. 

The optimal control problem associated with maximizing the net delivered useful mass can be formulated as,
\begin{equation}
\begin{array}{ll}
\underset{P_{\text{BL}}, \hat{\bm{\alpha}}, P_E}{\operatorname{minimize}} & J = -m_u, \\
\\
\text { subject to } & \dot{\boldsymbol{x}}(t)=\boldsymbol{f}(\boldsymbol{x}, \boldsymbol{a}), \\
& \boldsymbol{x}\left(t_0\right)=\boldsymbol{x}_0, \quad \boldsymbol{x}\left(t_f\right)=\boldsymbol{x}_f,  \quad m(t_0) = m_0,\\
& \|\bm{\hat{\alpha}}\| = 1, \\ 
& 0 \leq P_{\text{E}}(t) \leq P_\text{ava}(t,r), \\
& \boldsymbol{x}_\text{LB} \leq \boldsymbol{x} \leq \boldsymbol{x}_\text{UB}, \\
& m_\text{LB} \leq m \leq m_\text{UB}, \\
& P_\text{LB} \leq P_{\text{BL}} \leq P_\text{UB},
\end{array}
\label{eq:OCP}
\end{equation}
where the states and the mass are bounded between lower and upper bounds, denoted by subscripts, ``LB'' and ``UB'', respectively, and $P_{\text{E}}$ is bounded between a minimum value and the available power for the engine, $P_\text{ava}$. Calculation of the $P_\text{ava}$ value is explained in the next section. 

To size the solar arrays,  $P_{\text{BL}}$ is introduced as an additional optimization variable, which is the beginning-of-life power of the solar arrays. 
The optimization problem maximizes the net delivered useful mass, $m_u$, defined as,
\begin{equation}\label{eq:m_u}
    m_u = m_0 - \gamma_1P_{\text{BL}}-\gamma_2P_{\text{max}} - (1+\alpha_{\text{tk}})m_p
\end{equation}
where $m_0$ is the (fixed) initial spacecraft mass; $\gamma_1$ and $\gamma_2$ are the specific masses of solar arrays, and the power regulator unit and electric propulsion system (EPS), and $\alpha_{\text{tk}}$ is the tank mass ratio coefficient, $m_p$ is the fuel consumption calculated as $m_p = m_0 - m_f$. The cost can be further simplified by removing the constant terms without affecting the optimality. The optimal control problem, given in Eq.~\eqref{eq:OCP}, is converted into an NLP problem by discretizing the dynamics (adopting an RK4 transcription), states, and controls. The obtained NLP problem can be written as,

\begin{equation}
\begin{array}{ll}
\underset{P_{\text{BL}}, P_{E,i}, \hat{\bm{\alpha}}_i, \bm{X}_j}{\operatorname{minimize}} & J = -\left[m_1 - \gamma_1P_\text{BL} - \gamma_2P_\text{max} - (1+\alpha_\text{tk})(m_1 - m_N)\right], \\
\\
\text { subject to } & \bm{X}_{i+1} = \bm{X}_i + \frac{h_i}{6}(\bm{k}_1 + 2\bm{k}_2 + 2\bm{k}_3 + \bm{k}_4), \\
& \bm{x}_1=\boldsymbol{x}_0, \quad \boldsymbol{x}_N=\boldsymbol{x}_f,  \quad m_1 = m_0,\\
& \|\bm{\hat{\alpha}}_i\| = 1, \\
&0 \leq P_{\text{E},j} \leq P_\text{ava}(t_j,r_j), \\
& P_{\text{E},i} \leq P_\text{ava}(t_{i+1},r_{i+1}), \\
& \boldsymbol{x}_\text{LB} \leq \boldsymbol{x}_j \leq \boldsymbol{x}_\text{UB}, \\
& m_\text{LB} \leq m_j \leq m_\text{UB}, \\
& P_\text{LB} \leq P_{\text{BL}} \leq P_\text{UB},
\end{array}
\label{eq:direct_OCP}
\end{equation}
where, $\bm{X}_j = \bm{X}(t_j)$ with $\bm{X} = [\bm{x}^{\top}, m]^\top$,  $i = 1, \cdots, N-1$, $j = 1, \cdots, N$ with $N$ denoting the total number of nodes. The step size is defined as $h_i = t_{i+1}-t_{i}$. In Eq.~\eqref{eq:direct_OCP}, $\boldsymbol{x}_0$ and $\boldsymbol{x}_f$ denote the sets of initial (at the departure) and final (upon arrival) MEEs, respectively, and $r$ denotes the spacecraft distance from the Sun. 

It should be noted that since an NLP problem is formulated,  all the states are discretized and are available to form any cost function, which is of interest/importance. Therefore, the cost can be calculated using the final mass of the spacecraft from the discretized state, $m_N$, and the optimization problem can be written in a Mayer form (i.e., a terminal cost) instead of the Lagrange form that is used within the indirect optimization methods \cite{petukhov_joint_2019,arya_electric_2021}. 

\subsection{Power and Thruster Models}
The power generated by the solar arrays, $P_{\mathrm{SA}}$, is modeled as, \cite{arya_electric_2021, ellison_analytic_2018}
\begin{equation}
P_{\mathrm{SA}}=\frac{P_{\mathrm{BL}}}{r^2}\left[\frac{d_1+d_2 r^{-1}+d_3 r^{-2}}{1+d_4 r+d_5 r^2}\right](1-\sigma)^t,
\end{equation}
where $P_\text{BL}$ is the beginning-of-life power value of the solar arrays at 1 astronomical unit (au) distance from the Sun,  $\sigma \in [0.02,0.05]$ is a time-degradation factor (percent/year) and $t$ is the time passed from starting of the mission in years; $[d_1, d_2, d_3, d_4, d_5]$ are coefficients that are determined empirically, but in this work we consider $[1.1063, 0.1495, -0.299, -0.0432, 0]$. The power model considers the heliocentric distance as well as time degradation of the solar arrays~\cite{arya2021electric}. The available power for the thrusters is calculated as,
\begin{equation}
P_\text{ava}= \begin{cases}\eta_d P_{\max }, & \text { for } P_{\mathrm{SA}} \geq P_{\mathrm{sys}}+P_{\max }, \\ \eta_d\left(P_{\mathrm{SA}}-P_{\mathrm{sys}}\right), & \text { for } P_{\mathrm{SA}}<P_{\mathrm{sys}}+P_{\mathrm{max}},\end{cases}
\label{eq:powconst}
\end{equation}
where $P_\text{SA}$ is the power generated by solar arrays, $P_\text{sys}$ is the power required for other subsystems, $ P_{\text{max}}$ is the maximum power that the power processing unit can handle, and $\eta_d$ is the duty cycle of the thrusters. Due to the piece-wise continuous nature of the available power \cite{arya_electric_2021}, an L2 norm smoothing \cite{taheri2023l2} is used for a one-parameter smooth representation of the available power as,
\begin{equation}
\begin{aligned}
P_{\text{ava}} & =\eta_d\left[\chi_p P_{\mathrm{max}}+\left(1-\chi_p\right)\left(P_{\mathrm{SA}}-P_{\mathrm{sys}}\right)\right], \\
g_p &= P_\text{max} - (P_\text{sys} + P_\text{max}),\\
\chi_p & =\frac{1}{2}\left[1+\frac{g_p}{\sqrt{g_p^2 + \rho_p^2}}\right],
\end{aligned}
\end{equation}
where $\rho_p$ denotes a ``power'' smoothing parameter. We consider an SPT-140 Hall thruster~\cite{manzella1997performance}. This thruster throttle table lists 21 modes of operation, characterized by distinct thrust levels, mass flow rate, input power, $I_\text{sp}$, and efficiency parameters. The throttle table is given in Table~\ref{tab:spt140}.

\begin{table}[!htbp]
\caption{Throttle table for SPT-140 Hall thruster~\cite{manzella1997performance}.}
    \centering
    \begin{tabular}{cllllc}
\hline $\begin{array}{c}\text { Mode } \\
\#\end{array}$ & $\begin{array}{l}\text { Power } \\
\text { (Watts) }\end{array}$ & $\begin{array}{l}T \\
(\mathrm{mN})\end{array}$ & $\begin{array}{l}\dot{m} \\
(\mathrm{mg} / \mathrm{s})\end{array}$ & $\begin{array}{l}I_{\mathrm{sp}} \\
(\mathrm{sec})\end{array}$ & $\eta$ \\
\hline 1 & 4989 & 263 & 13.9 & 1929 & 0.5 \\
2 & 4620 & 270 & 16.5 & 1670 & 0.48 \\
3 & 4589 & 287 & 17.8 & 1647 & 0.5 \\
4 & 4561 & 264 & 16.4 & 1645 & 0.47 \\
5 & 4502 & 260 & 16.2 & 1641 & 0.46 \\
6 & 4375 & 246 & 14 & 1790 & 0.49 \\
7 & 3937 & 251 & 17.5 & 1461 & 0.46 \\
8 & 3894 & 251 & 17.5 & 1464 & 0.46 \\
9 & 3850 & 251 & 17.5 & 1464 & 0.47 \\
10 & 3758 & 217 & 13.9 & 1597 & 0.45 \\
11 & 3752 & 221 & 13.9 & 1617 & 0.47 \\
12 & 3750 & 215 & 13.6 & 1614 & 0.45 \\
13 & 3460 & 184 & 17.1 & 1099 & 0.29 \\
14 & 3446 & 185 & 20.4 & 925 & 0.24 \\
15 & 3402 & 189 & 16.3 & 1181 & 0.32 \\
16 & 3377 & 201 & 15.8 & 1302 & 0.38 \\
17 & 3376 & 175 & 18.2 & 979 & 0.25 \\
18 & 3360 & 198 & 14.7 & 1371 & 0.4 \\
19 & 3142 & 191 & 13.8 & 1409 & 0.42 \\
20 & 3008 & 177 & 11.4 & 1579 & 0.46 \\
21 & 1514 & 87 & 6.1 & 1449 & 0.41 \\
\hline
\end{tabular}
    
    \label{tab:spt140}
\end{table}

\subsection{Smooth Mode Selection}
To select the current mode, $\bm P_\text{E} = [P_{E,1},\cdots,P_{E,j}]$ is used as the engine power control, which represents the collection of input power to the engine at each time step. The engine power is bounded by the available power, $\bm P_{\text{ava}}$ and the available power is bounded by the maximum power that PPU can process, $\bm P_{\text{ava}} \leq P_{\text{max}}$. This way, the NLP solver can control $\bm P_\text{E}$ such that the selected mode maximizes the useful mass and the $\bm P_\text{E}$ satisfies the power constraints. The discontinuous mode selection is regularized by the L2-norm smoothing~\cite{taheri2023l2}. 

At each time instant, the $P_\text{E}(t_i)$ is compared against the list of selected mode powers, $\bm{P}_{\text{sel}} = [P_{\text{sel},1},\cdots, P_{\text{sel},N_\text{mode}}]^{\top} \in \mathbb{R}^{{N_\text{mode}}} $,  to determine the activated power-feasible mode. We drop the argument for the sake of simplicity and refer to the value of $P_\text{E}(t_i)$ as $P_\text{E}$. 
The comparison of the $P_\text{E}$ with the first element of $\bm{P}_{\text{sel}}$, i.e., $P_{\text{sel},1}$, is achieved by following the steps outlined in Eq.~\eqref{eq:psel1comp} as,
\begin{equation} \label{eq:psel1comp}
\begin{aligned}
g_{1,1}&=P_\text{E}- P_{\text {sel}, 1}, \\
\zeta_{1,1}&=\frac{1}{2}\left(1+\frac{g_{1,1}}{\sqrt{g_{1,1}^2 + \rho_{e}^2}}\right), \\
\eta_{1} &= \zeta_{1,1},
\end{aligned}
\end{equation}
where $g_{1,1}$ is the power-sufficiency switching function for mode 1, $\rho_e$ is its associated smoothing parameter, and $\eta_{1,1}$ is the first element of an activation weight vector. If $P_\text{E} > P_{\text{sel},1}$, the first mode gets activated, $\zeta_{1,1} = \eta_{1} = 1$ as $\rho_e \rightarrow 0$. The remaining mode powers should be compared against the $P_\text{E}$ as well. The previous and current power values are important for the remaining modes to avoid the simultaneous activation of multiple modes. Two power switching functions are defined and the remaining elements of the activation weight vector are obtained by following the steps, 
\begin{equation}
\begin{aligned}
g_{i, 1}&=P_\text{E}-P_{\text {sel},i-1}, \\
g_{i, 2}&=P_\text{E}-P_{\text {sel}, i}, \\
\zeta_{i, 1}&=1-\frac{1}{2}\left(1+\frac{g_{i,1}}{\sqrt{g_{i,1}^2 + \rho_{e}^2}}\right), \\
\zeta_{i, 2}&=\frac{1}{2}\left(1+\frac{g_{i,2}}{\sqrt{g_{i,2}^2 + \rho_{e}^2}}\right),\\
\eta_{i} &= \zeta_{i, 1}\zeta_{i, 2}, \\
\end{aligned}
\end{equation}
where $i = 2, \cdots, N_\text{mode}$. Please note that $\zeta_{i, 1}$ determines if the previous mode power value is higher than the $P_\text{E}$. If that is the case, 
current mode is not activated; otherwise,  the current power value is still compared against $P_\text{E}$ to determine whether it can be activated or not. 

For instance, assume $\eta_1 = 1$ with $P_\text{E} = 4000$ W, $P_{\text {sel},1} = 3500$ W, and $P_{\text {sel},2} = 3000$ W. Even though both $g_{2,1}$ and $g_{2,2}$ are positive, $\zeta_{2,1} = 0$ and $\zeta_{2,2} = 1$. Therefore, 
$\eta_2 = 0$ since both $P_\text{E} > P_{\text {sel},1}$ and $P_\text{E} > P_{\text {sel},2}$. This logic is configured to ensure that only one of the elements of the $\bm{\eta}$ vector becomes 1 (in the limit as $\rho_e \rightarrow 0$). Using the activation vector elements, we construct the activation weight vector such that $\bm{\eta} = [\eta_1, \eta_2, \cdots \eta_{N_\text{mode}}]^\top \in \mathbb{R}^{N_\text{mode}}$. The activation weight vector is used to determine the thrust and mass flow rate values, which appear in Eqs.~\eqref{eq:dyn} and ~\eqref{eq:mdot}, 
\begin{equation}
\begin{aligned}
T&=\bm{\eta}^\top \bm{T}_{\text {sel}}, \quad & 
\dot{m}&=-\bm{\eta}^\top \dot{\bm{m}}_{\text {sel}},
\end{aligned}
\end{equation}
where $\bm{T}_{\text {sel}}$ and $\dot{\bm{m}}_{\text {sel}}$ denote the sub-set of modes from the third and fourth columns of Table~\ref{tab:spt140}. The current mode, $n$, as defined in Eq.~\eqref{eq:fmode}, is the element of $\bm{\eta} = 1$, when $\rho_e \rightarrow 0$. 

\subsection{Initial Guess Generation and Continuation}
The initial guess is generated at the very first step of a continuation procedure to reduce the values of the smoothing parameters, $\rho_p$ and $\rho_e$. In the subsequent steps of the numerical continuation procedure, the previous solution is used as an initial guess. For the very first step, the initial guess is needed for $\bm{x}$, $m$, $\bm{u}$, and $P_\text{BL}$. For the states, $\bm x$ and $m$, linear interpolation is performed for the initial guess generation. For $\bm u = [\hat{\bm{\alpha}}, P_\text{E}]$, there are two steps involved. The thrust steering direction ($\hat{\bm{\alpha}}$) is chosen parallel to the velocity direction of the interpolated states. For the engine power control, an initial guess, $P_\text{E}$, $P_\text{ava}$ is calculated with the interpolated steps, and after an estimate for the distance from the sun is obtained. Then, $P_\text{ava}$ is used as an input to the smooth mode-selection algorithm to generate $\bm \eta$ to calculate $\bm{\eta}^\top \bm{P}_\text{sel}$. 

The smooth mode-selection algorithm switching functions are scaled with the highest power of the $\bm{P}_\text{sel}$ to regularize the values of the switching functions to improve convergence. The continuation strategy is to lower both $\rho_p$ and $\rho_e$ simultaneously. If a large step is taken towards decreasing $\rho_e$, we keep the value of $\rho_p$ to be constant. Another strategy that we considered is to keep the steps small such that the solver quickly converges without performing a significant number of iterations. If the solver takes a long time to converge, it is a sign that the considered step size is too large. In our computational experience, setting the number of discretization points to 100 or 200 nodes makes the solver converge faster for the considered problems.

\section{Numerical Results}
The Earth to Comet-67P time-fixed rendezvous mission is considered. We consider four sub-problems corresponding to one to four SPT-140 thruster modes. In addition, an extra 22nd mode with $0$ N thrust level and zero value of propellant consumption rate is considered to represent coast arcs. An NLP problem is formulated using the CasADi tool on MATLAB~\cite{andersson2019casadi}. The formulated problem is solved using IPOPT solver~\cite{biegler2009large}. The fourth-order Runge-Kutta method is used for discretizing the dynamics. The number of nodes used in the solutions is $N = 300$. 

Upon using the canonical units, the boundary conditions in MEE states are as follows corresponding to 2 revolutions of the spacecraft around the Sun:
\begin{equation*}
\begin{split}
\bm{x}_0 &= [0.998874284410563, -0.00294251935124146, 0.0164376759007608, \\
         &\quad -5.51480481733780e-06, 7.12277764431642e-06, 10.9784865869657], \\
\bm{x}_f &= [2.04295724237197, 0.292069230030979, 0.570126626743441, \\
         &\quad 0.0394123086323580, 0.0472705619148424, 28.3786463271836].
\end{split}
\end{equation*}

The initial mass of the spacecraft is fixed at $m_0 = 3000$ kg; the upper and lower bounds for $P_\text{BL} \in [10, 30]$ kW; $\sigma = 0.02$, $P_\text{max} = 4863$ W, $P_\text{sys} = 590$ W, $\gamma_1 = 10 \times 10^{-3}$ $\text{W}/\text{kg}$, $\gamma_2 = 15 \times 10^{-3}$ $\text{W}/\text{kg}$, $\eta_d = 0.95$. The fixed mission flight time is $1770$ days.

 The numerical results for the four cases are summarized in Table~\ref{tab:results}. The definitions of various mass expressions summarized in Table~\ref{tab:results} are given herein. The useful mass is given in Eq.~\eqref{eq:m_u} in which $m_\text{SA} = \gamma_1P_\text{BL}$ is the mass of the solar arrays, $m_\text{PSPU} = m_\text{SA} + \gamma_2P_\text{max}$ is the mass of the power supply and propulsion unit, and $m_\text{PSFS} = (1+\alpha_\text{tk})m_p$ is the mass of the propellant storage and feeding system. 
 As the considered number of operational modes increases, the useful mass capacity of the spacecraft increases, except in the three-mode case. Even though the useful mass of the three-mode solution is slightly smaller than the two-mode case, it is attributed to the ``unrealistic'' peaks in the $\bm P_\text{E}$ profile. 
 
 As the number of modes increases, the solar arrays get smaller, since $P_\text{BL}$ gets smaller. Therefore, considering more modes helped to design a spacecraft with a smaller solar array with a higher net delivered useful mass capacity. Compared to the one-mode case, fuel consumption decreases as more modes are considered. Having more modes gives more flexibility where a thrust arc does not necessarily correspond to the highest thrust level. Therefore, less fuel is consumed. The fuel consumption is almost the same for the three- and four-mode cases. However, the three-mode case results in heavier solar arrays. 
 
\begin{figure}[!htbp]
\centering
\begin{subfigure}{0.5\textwidth}
  \centering
\includegraphics[width=1\linewidth]{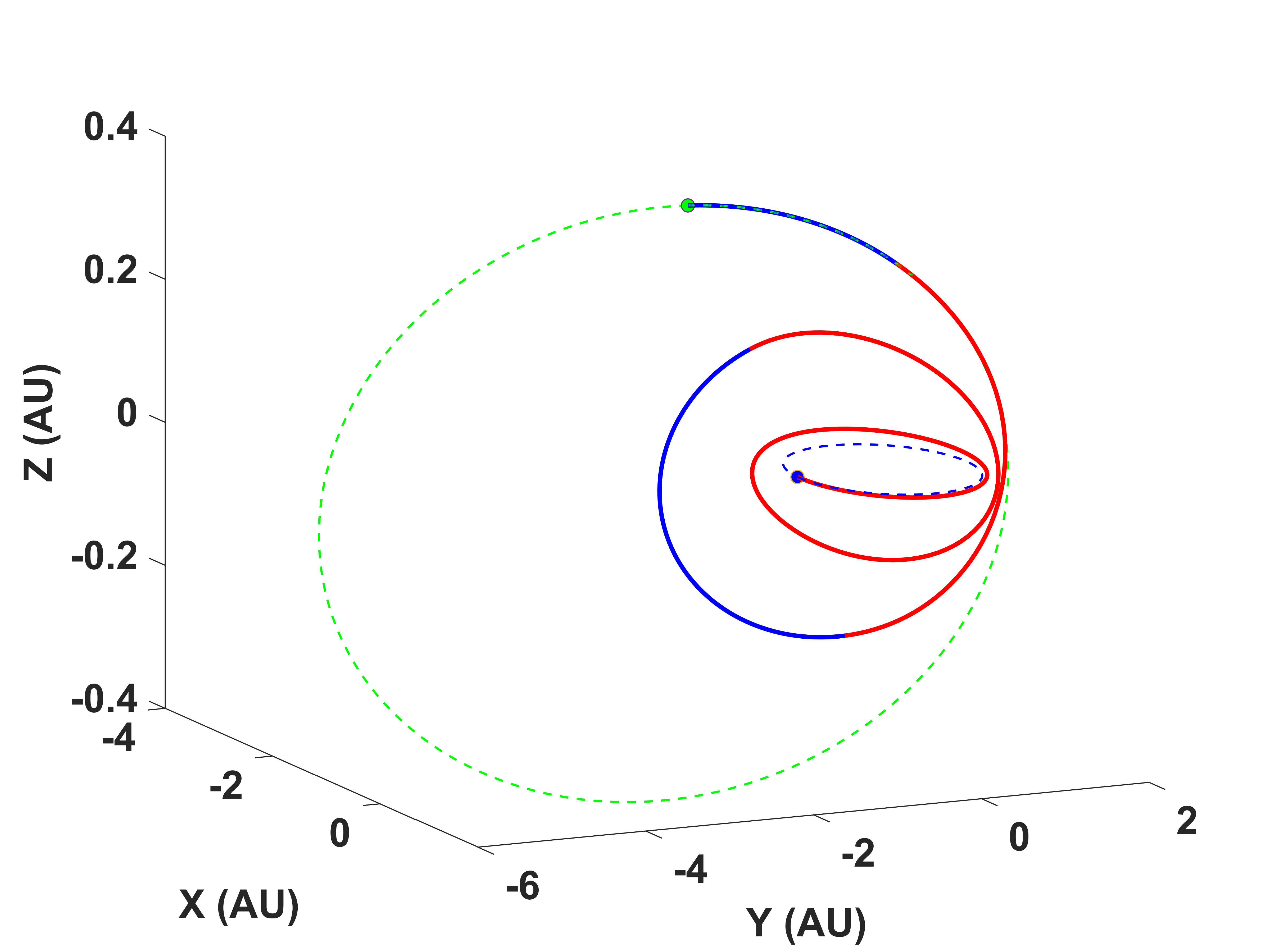}
  \caption{}
  \label{fig:traj_1mode}
\end{subfigure}%
\begin{subfigure}{0.5\textwidth}
  \centering
\includegraphics[width=1\linewidth]{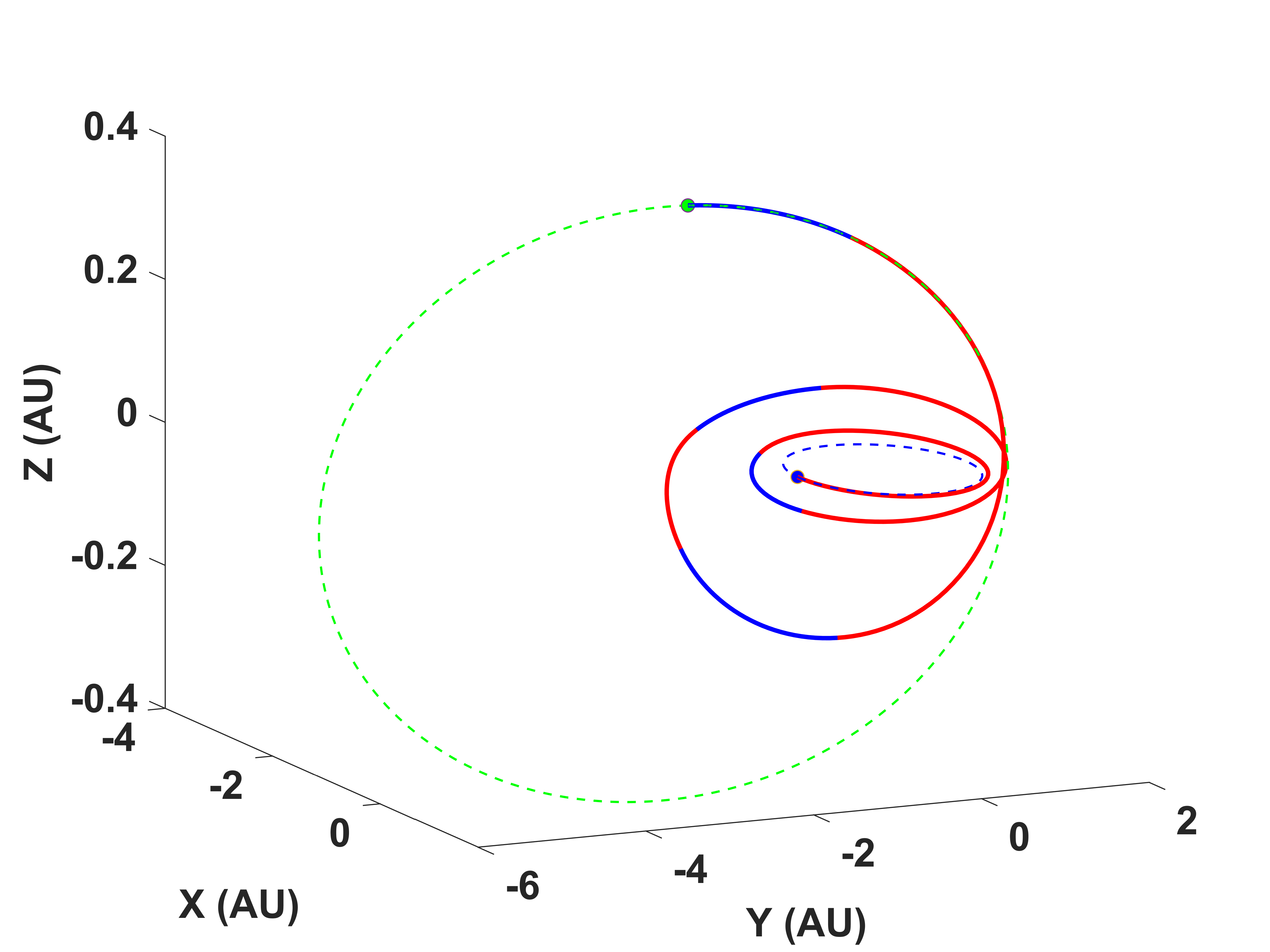}
  \caption{}
  \label{fig:traj_2mode}
\end{subfigure}
\begin{subfigure}{0.5\textwidth}
  \centering
\includegraphics[width=1\linewidth]{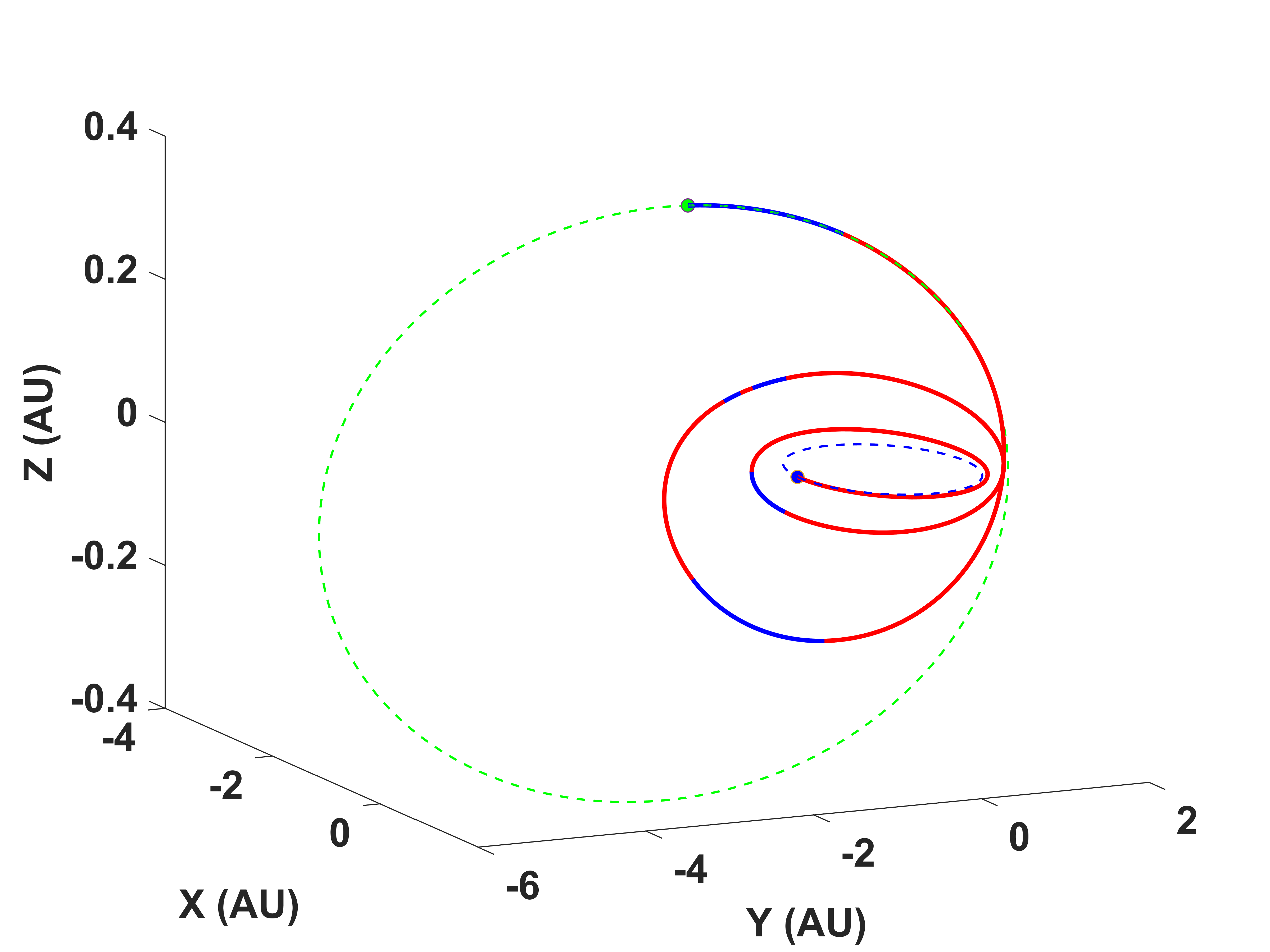}
  \caption{}
  \label{fig:traj_3mode}
\end{subfigure}%
\begin{subfigure}{0.5\textwidth}
  \centering
\includegraphics[width=1\linewidth]{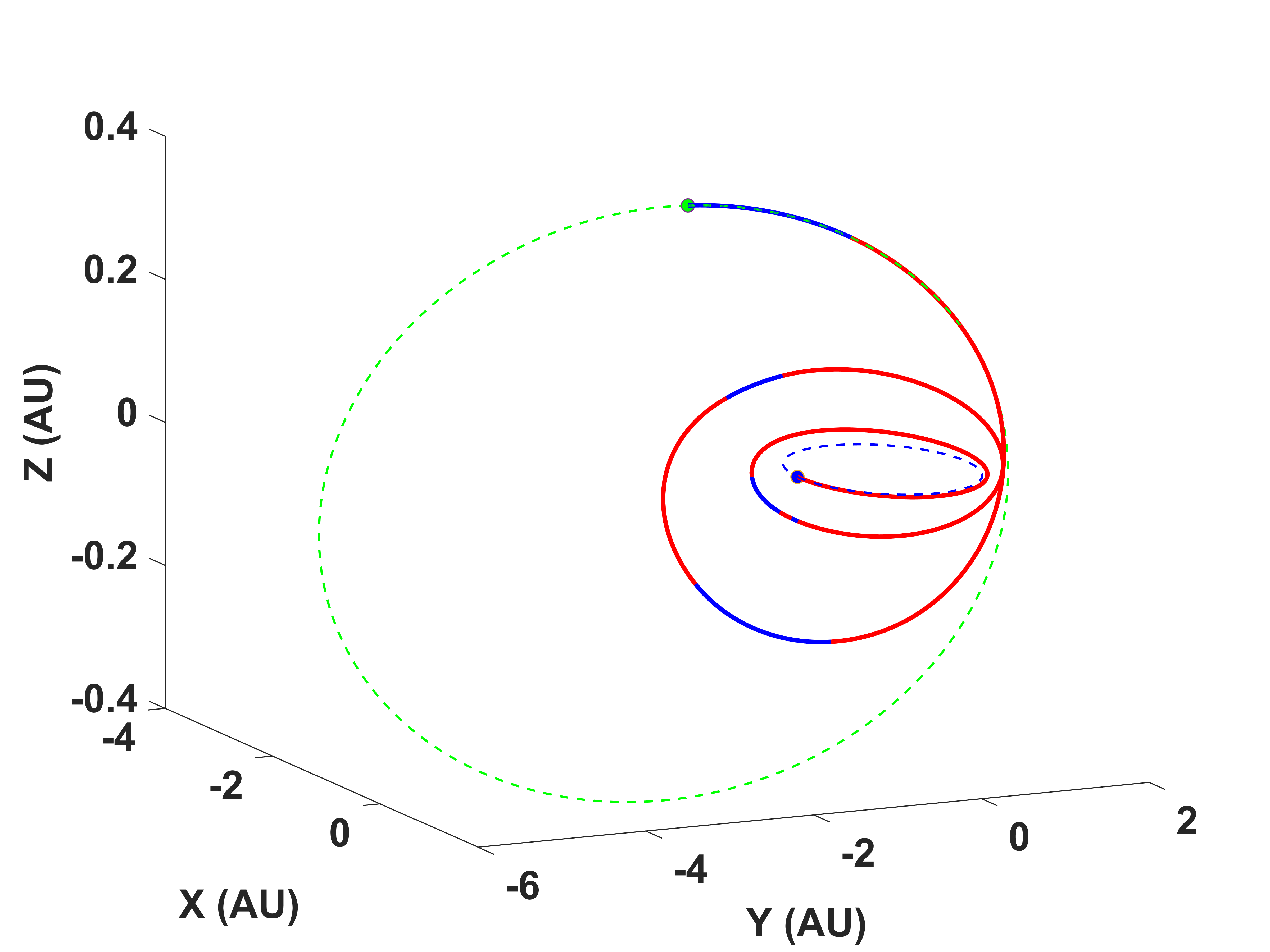}
  \caption{}
  \label{fig:traj_4mode}
\end{subfigure}
\caption{Earth-67P problem trajectories: (a) one mode, (b) two modes, (c) three modes, (d) four modes. }
\label{fig:traj}
\end{figure}

\begin{figure}[!htbp]
\centering
\begin{subfigure}{0.5\textwidth}
  \centering
\includegraphics[width=1\linewidth]{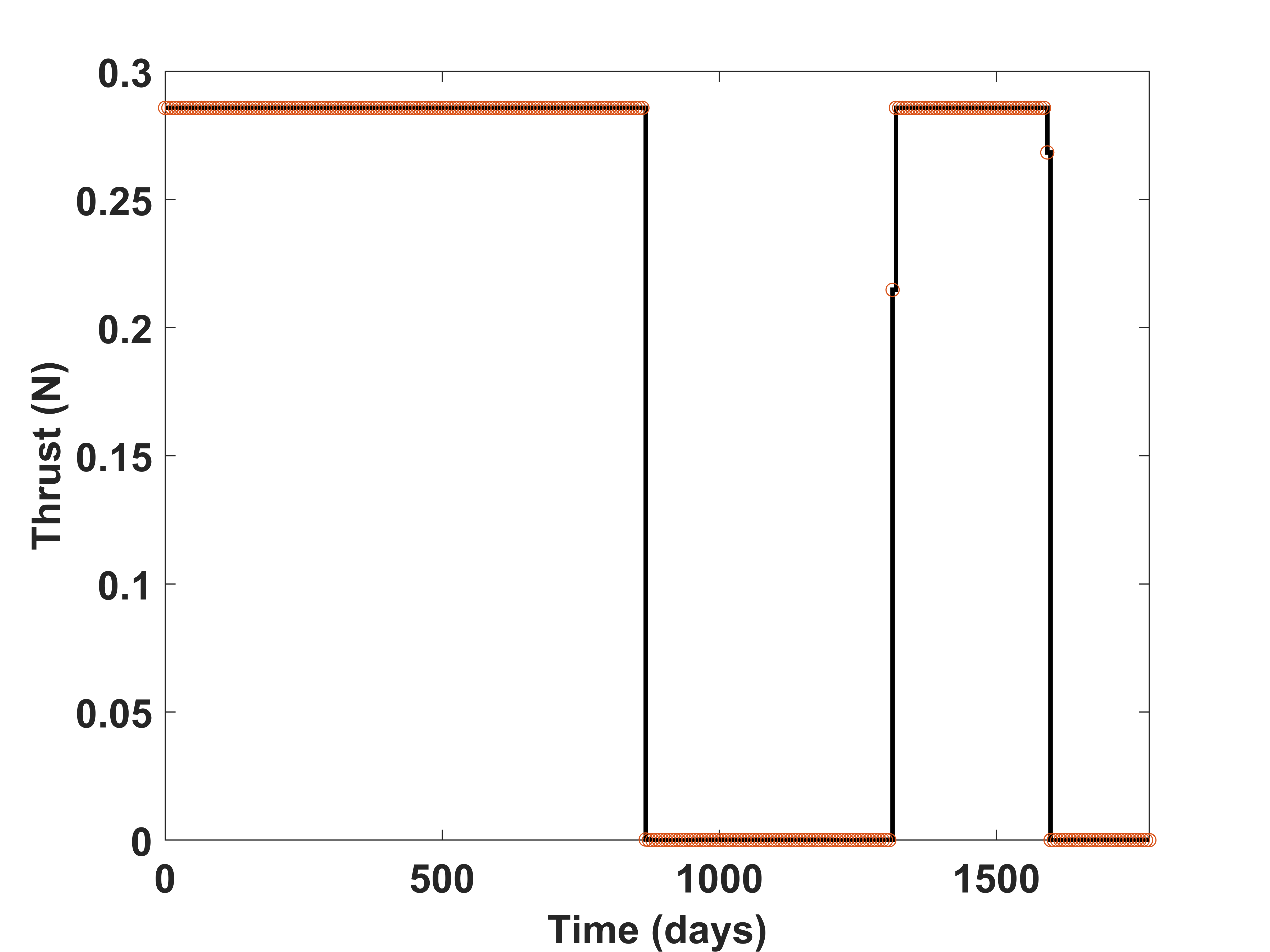}
  \caption{}
  \label{fig:thrust_1mode}
\end{subfigure}%
\begin{subfigure}{0.5\textwidth}
  \centering
\includegraphics[width=1\linewidth]{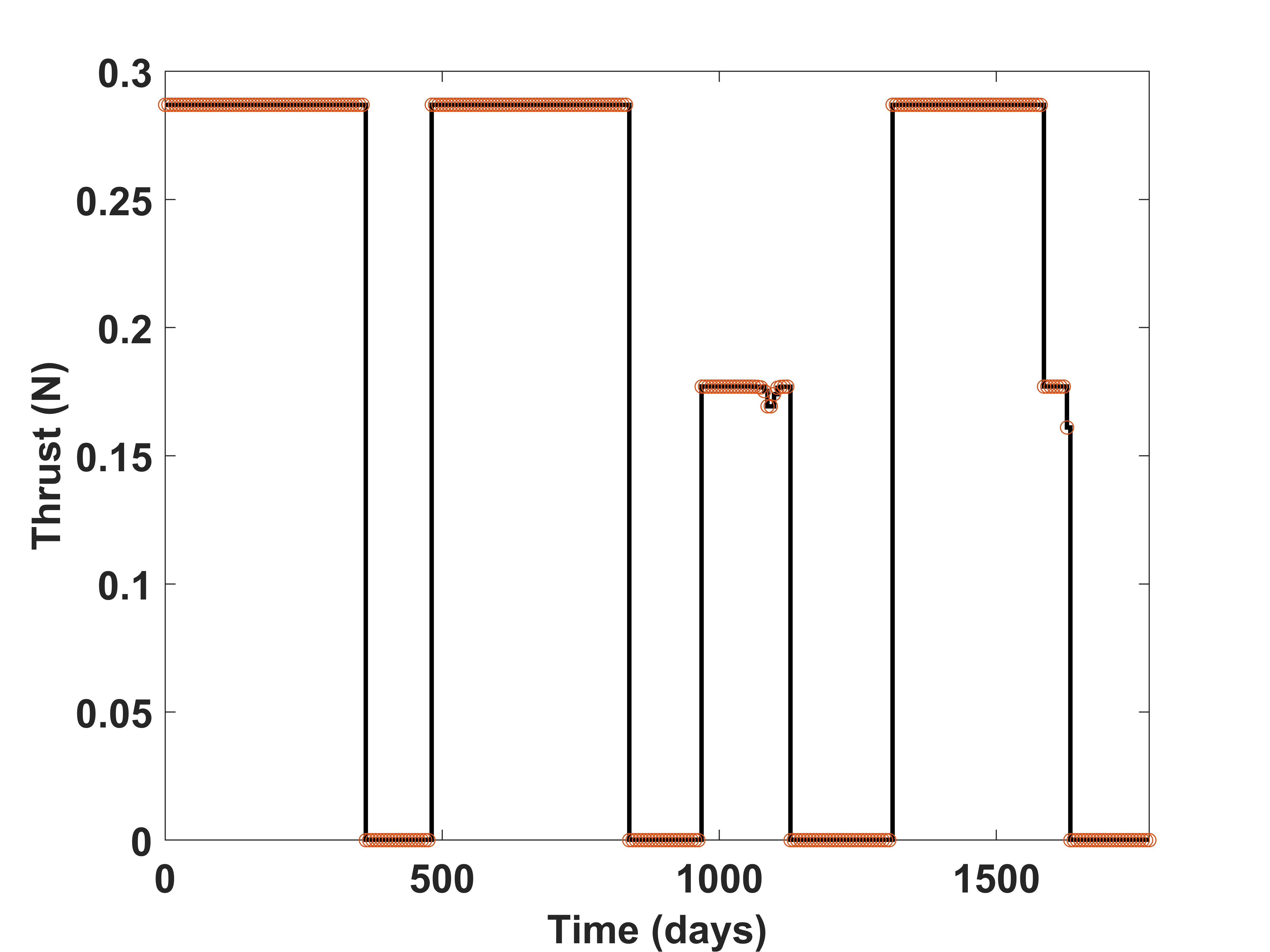}
  \caption{}
  \label{fig:thrust_2mode}
\end{subfigure}
\begin{subfigure}{0.5\textwidth}
  \centering
\includegraphics[width=1\linewidth]{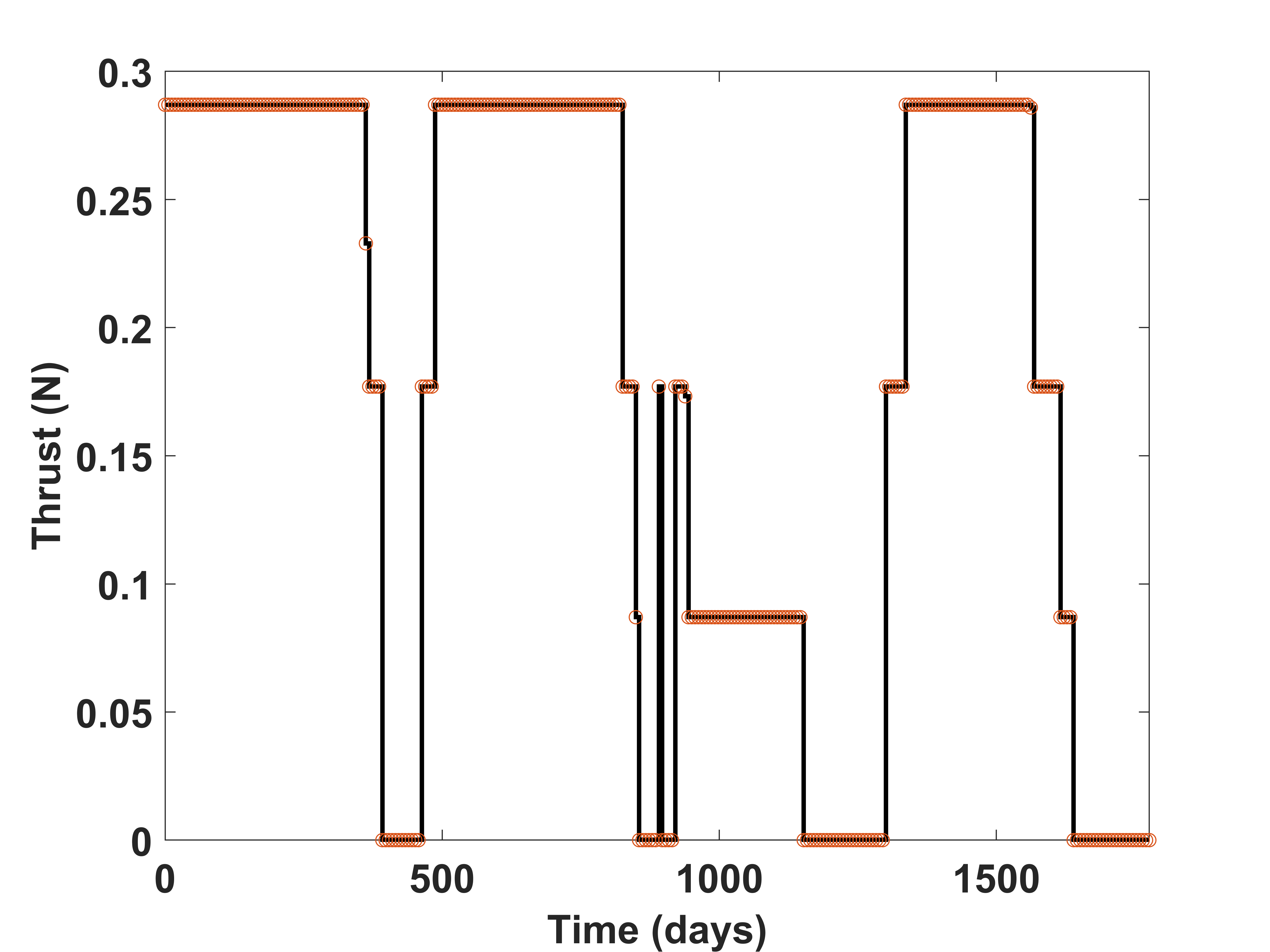}
  \caption{}
  \label{fig:thrust_3mode}
\end{subfigure}%
\begin{subfigure}{0.5\textwidth}
  \centering
\includegraphics[width=1\linewidth]{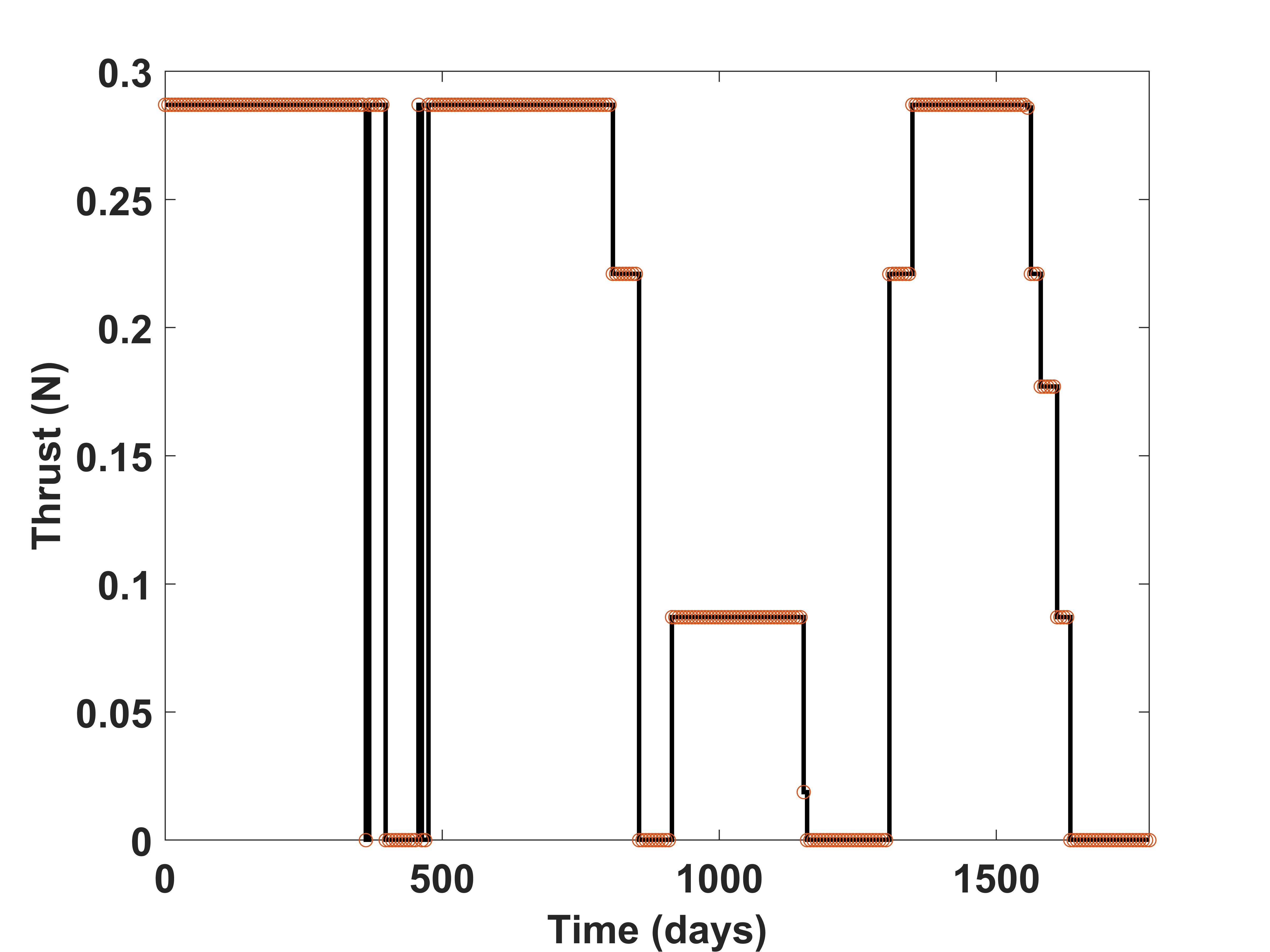}
  \caption{}
  \label{fig:thrust_4mode}
\end{subfigure}

\caption{Earth-67P problem thrust profiles: (a) one mode, (b) two modes, (c) three modes, and (d) four modes.}
\label{fig:thrust}
\end{figure}

\begin{figure}[!htbp]
\centering
\begin{subfigure}{0.5\textwidth}
  \centering
\includegraphics[width=1\linewidth]{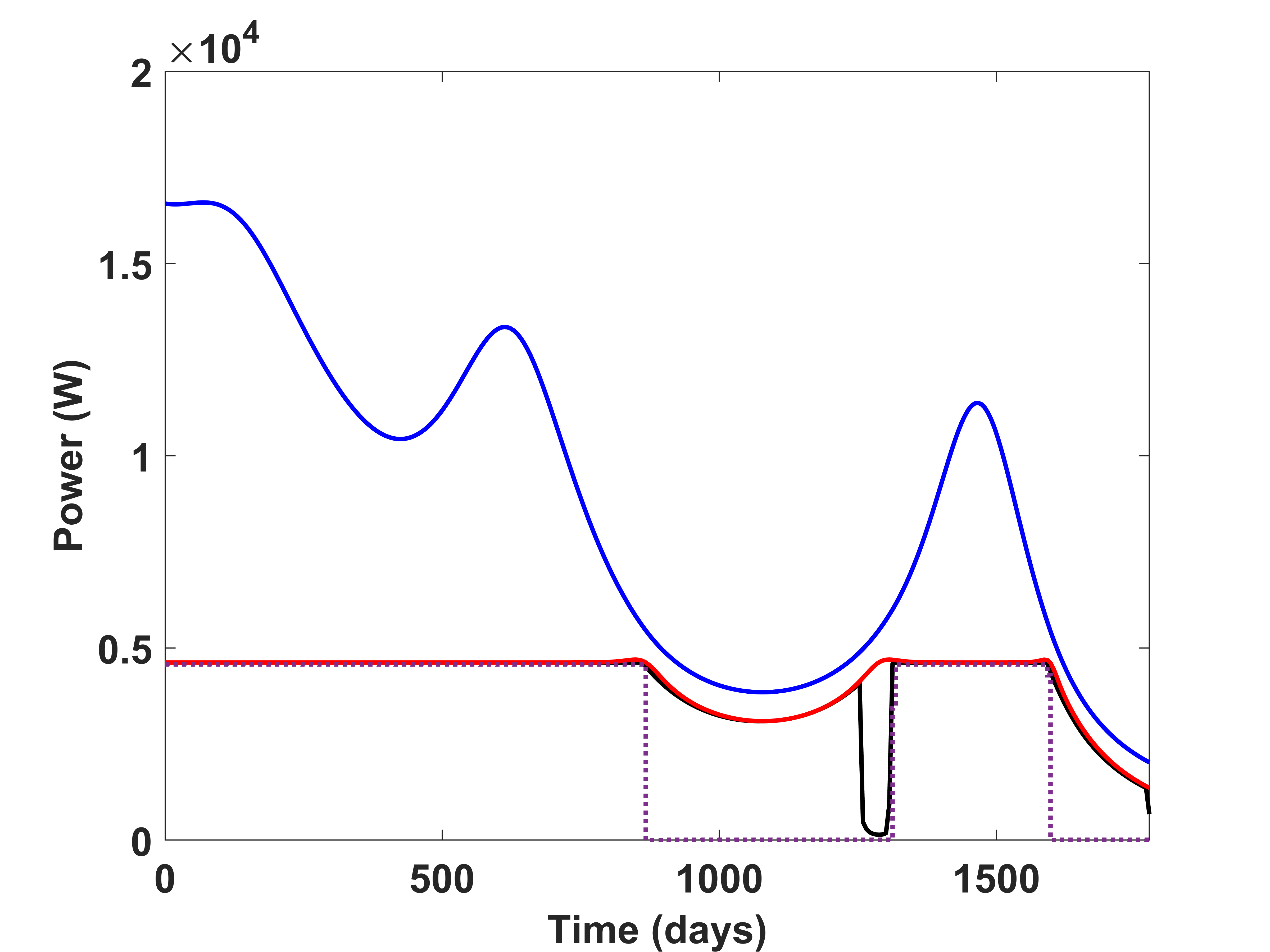}
  \caption{}
  \label{fig:power_1mode}
\end{subfigure}%
\begin{subfigure}{0.5\textwidth}
  \centering
\includegraphics[width=1\linewidth]{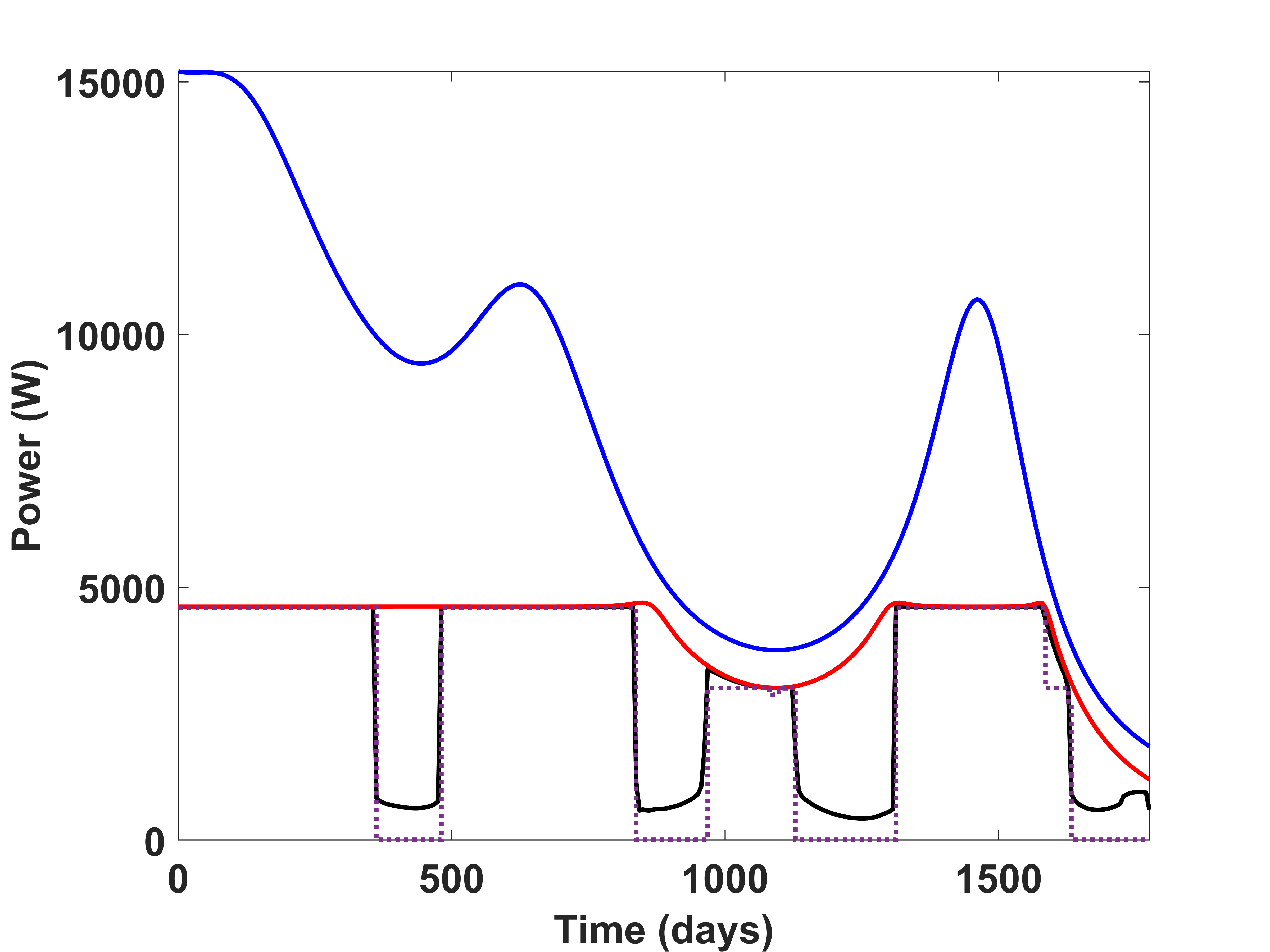}
  \caption{}
  \label{fig:power_2mode}
\end{subfigure}
\begin{subfigure}{0.5\textwidth}
  \centering
\includegraphics[width=1\linewidth]{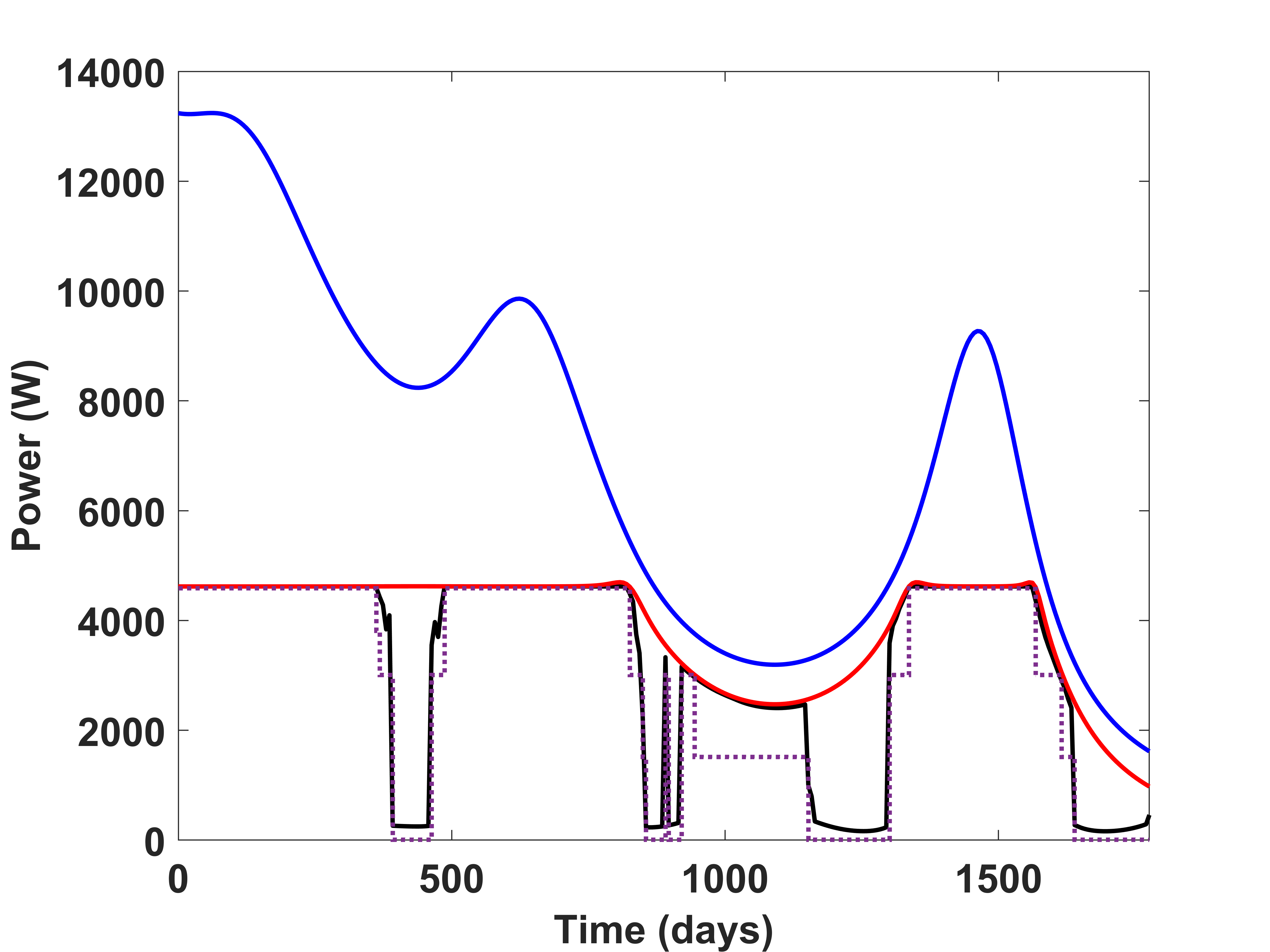}
  \caption{}
  \label{fig:power_3mode}
\end{subfigure}%
\begin{subfigure}{0.5\textwidth}
  \centering
\includegraphics[width=1\linewidth]{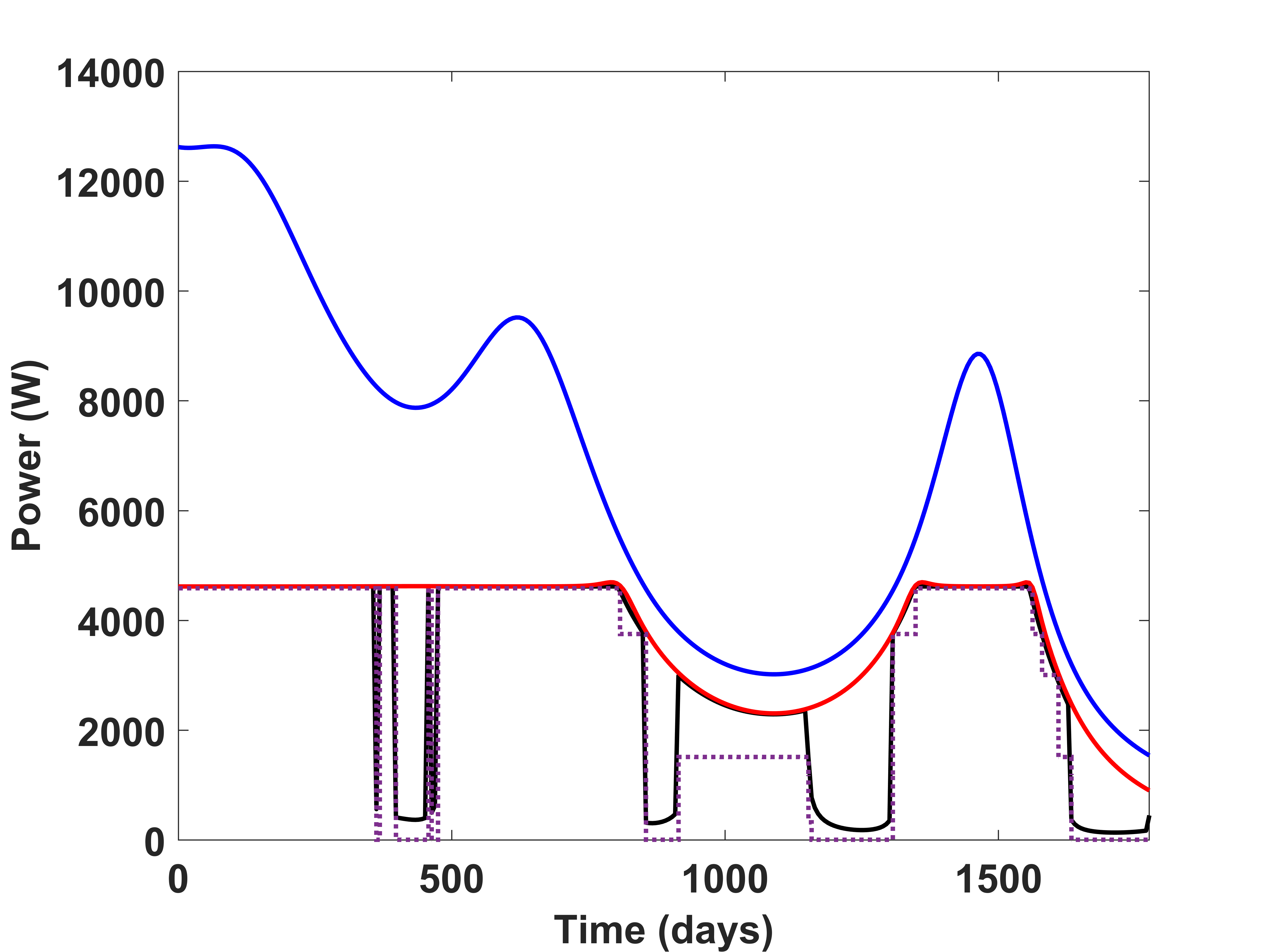}
  \caption{}
  \label{fig:power_4mode}
\end{subfigure}
\caption{Earth-67P problem power profiles: (a) one mode, (b) two modes, (c) three modes, and (d) four modes. The blue line is power generated by solar arrays, red is the available power, purple dashed is the selected power input and black is the $\bm P_\text{E}$.}
\label{fig:power}
\end{figure}

\begin{figure}[!htbp]
\centering
\begin{subfigure}{0.5\textwidth}
  \centering
\includegraphics[width=1\linewidth]{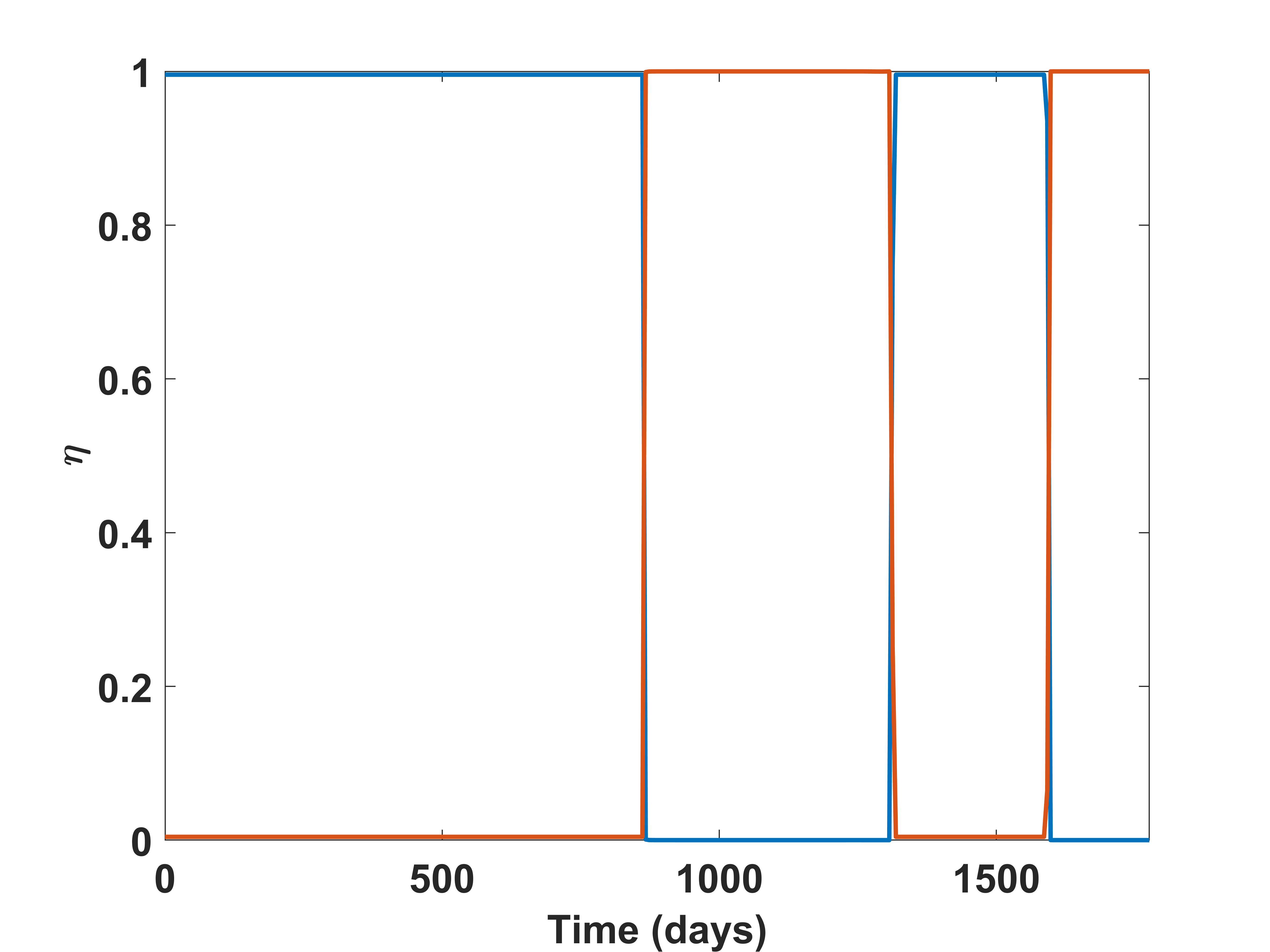}
  \caption{}
  \label{fig:eta_1mode}
\end{subfigure}%
\begin{subfigure}{0.5\textwidth}
  \centering
\includegraphics[width=1\linewidth]{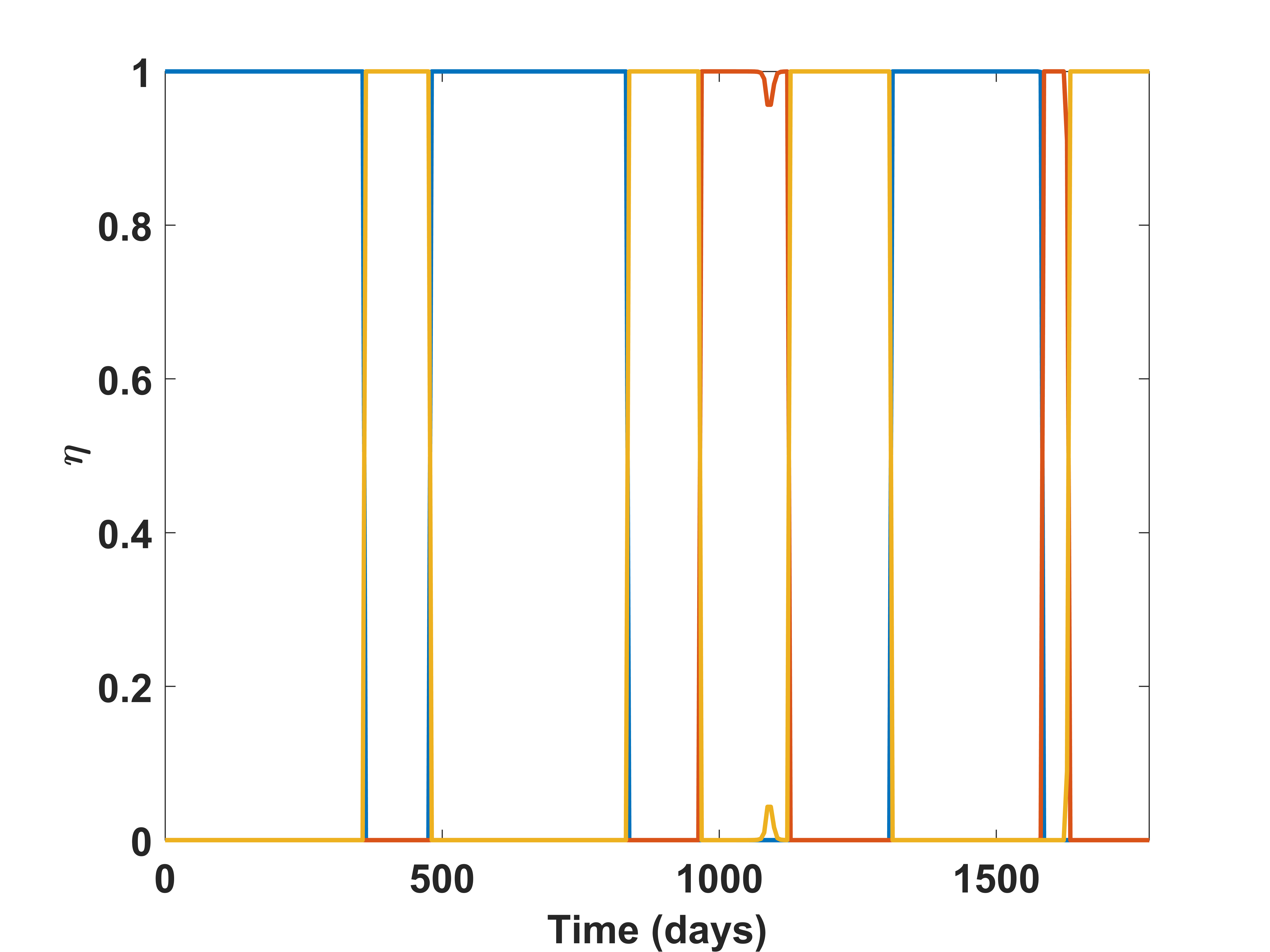}
  \caption{}
  \label{fig:eta_2mode}
\end{subfigure}
\begin{subfigure}{0.5\textwidth}
  \centering
\includegraphics[width=1\linewidth]{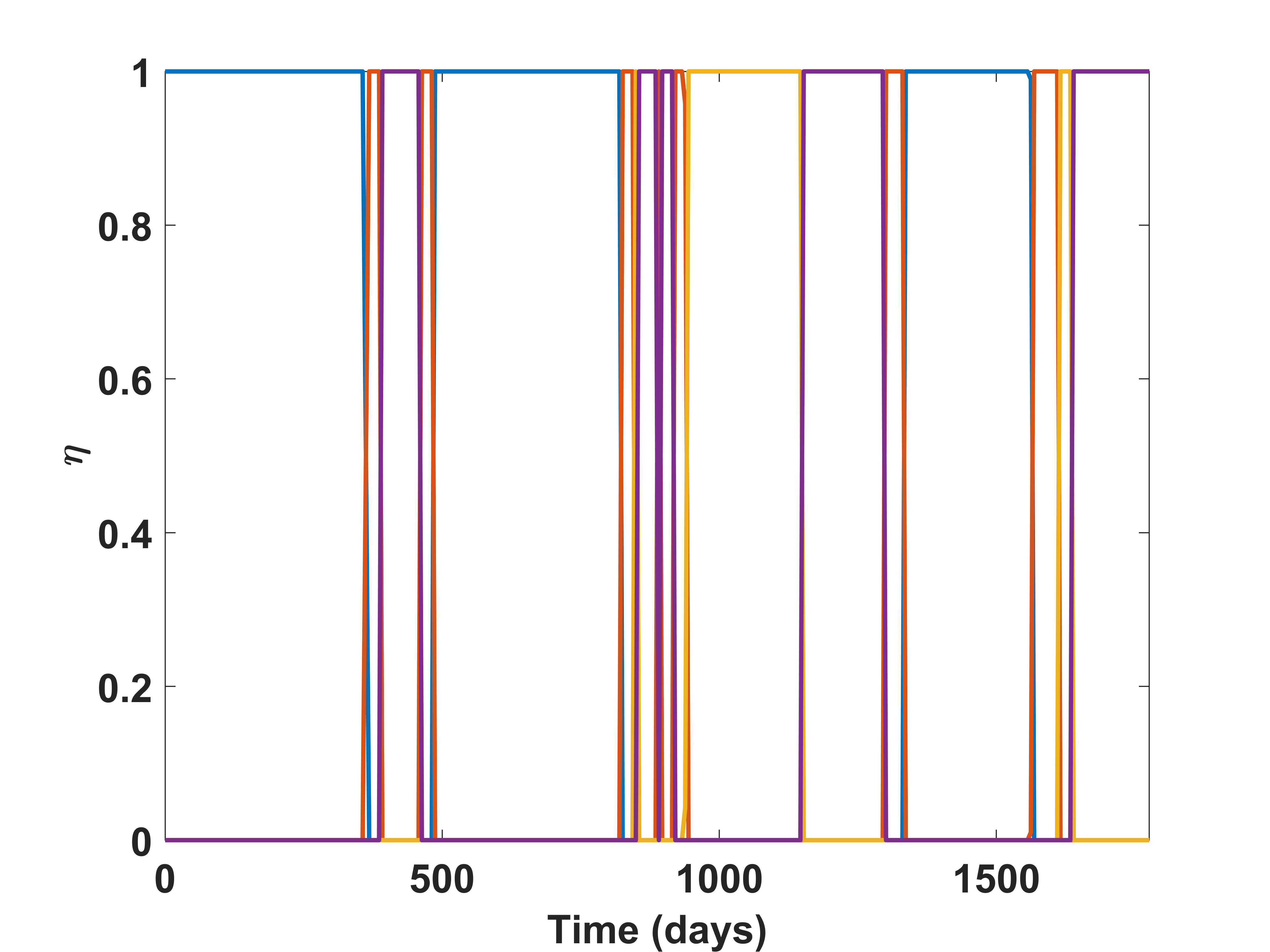}
  \caption{}
  \label{fig:eta_3mode}
\end{subfigure}%
\begin{subfigure}{0.5\textwidth}
  \centering
\includegraphics[width=1\linewidth]{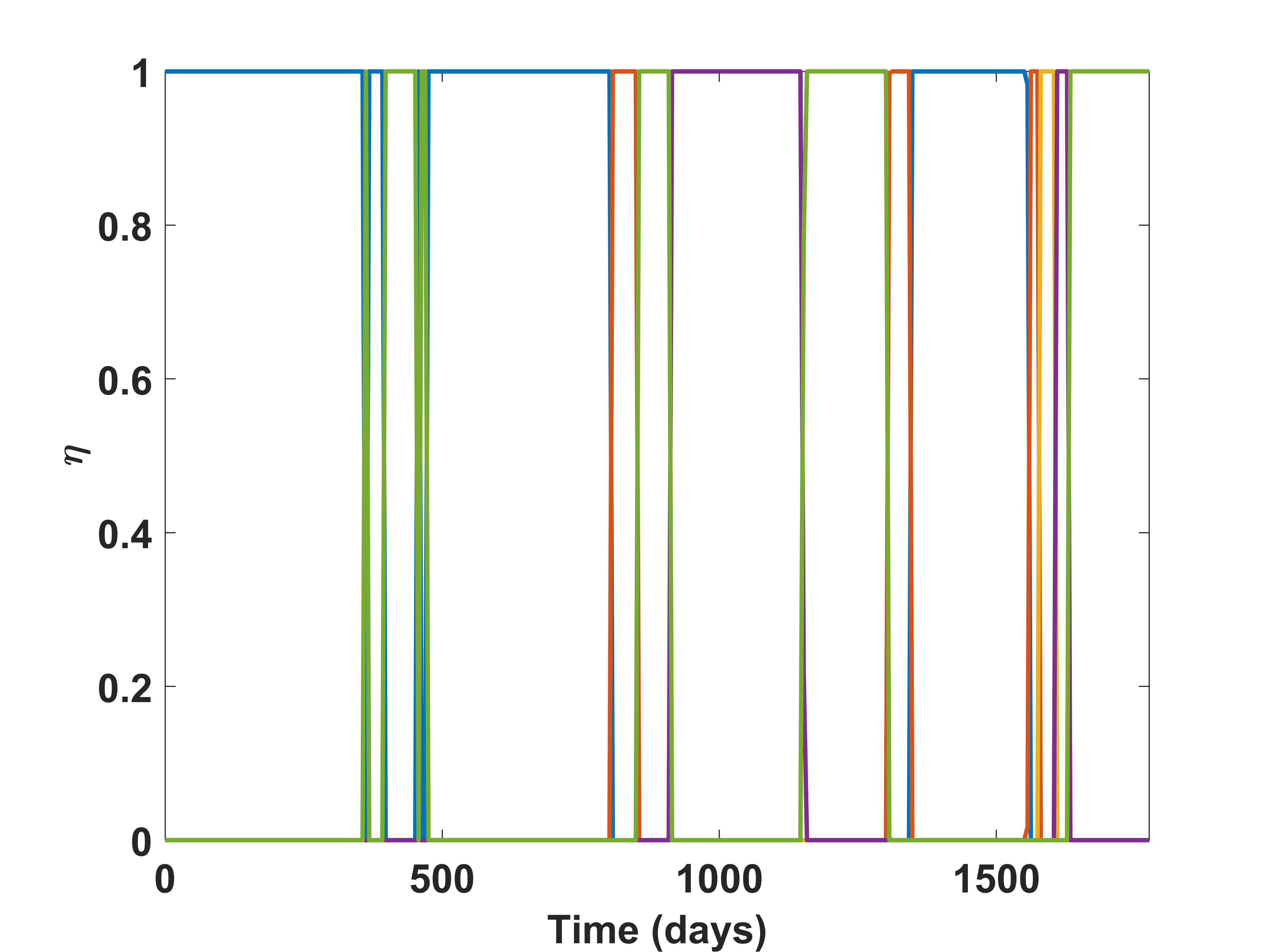}
  \caption{}
  \label{fig:eta_4mode}
\end{subfigure}

\caption{Earth-67P problem activation weight vectors ($\bm \eta$): (a) one mode, (b) two modes, (c) three modes, (d) four modes.}
\label{fig:eta}
\end{figure}

\begin{figure}[!htbp]
\centering
\begin{subfigure}{0.5\textwidth}
  \centering
\includegraphics[width=1\linewidth]{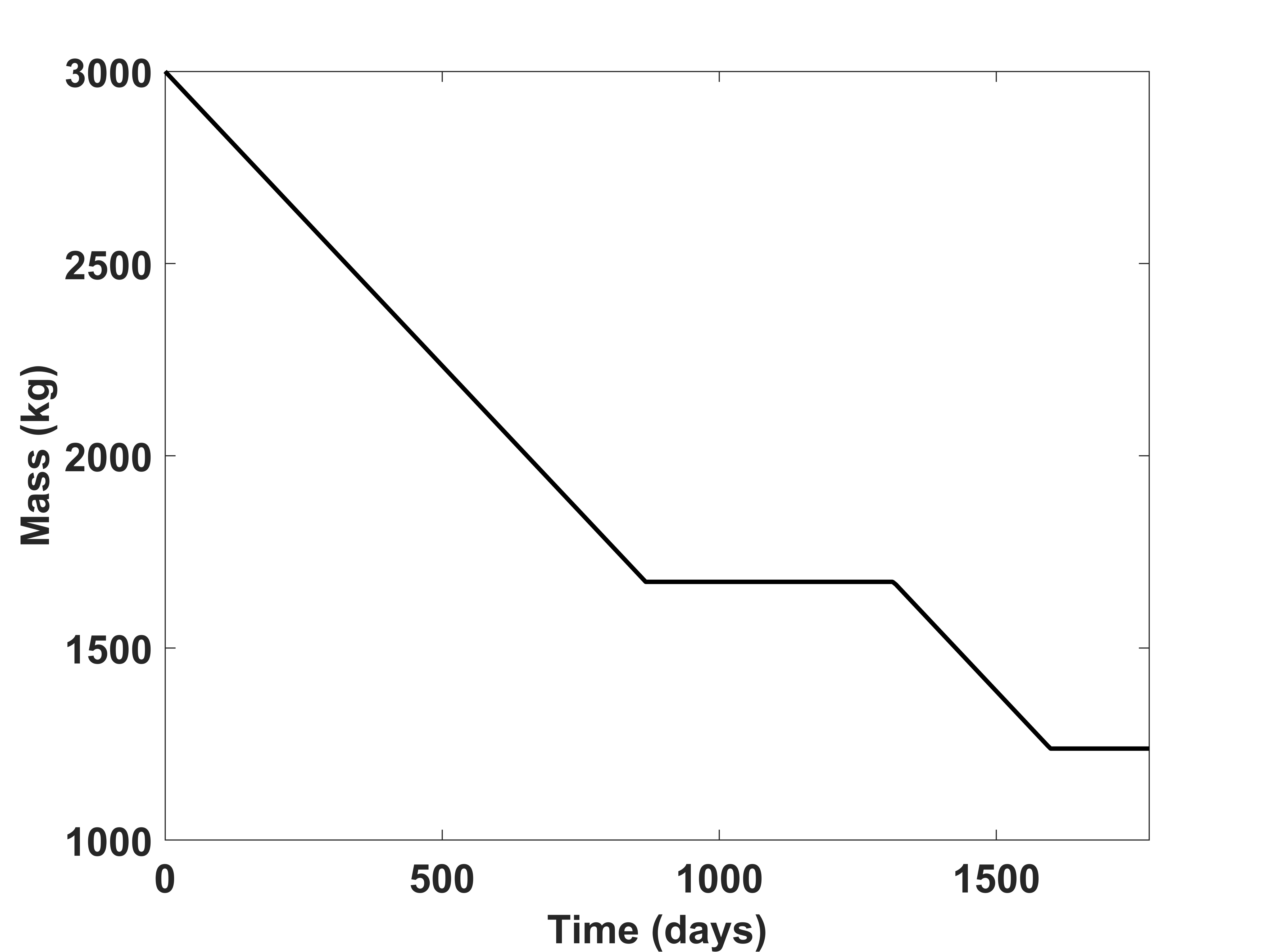}
  \caption{}
  \label{fig:mass_1mode}
\end{subfigure}%
\begin{subfigure}{0.5\textwidth}
  \centering
\includegraphics[width=1\linewidth]{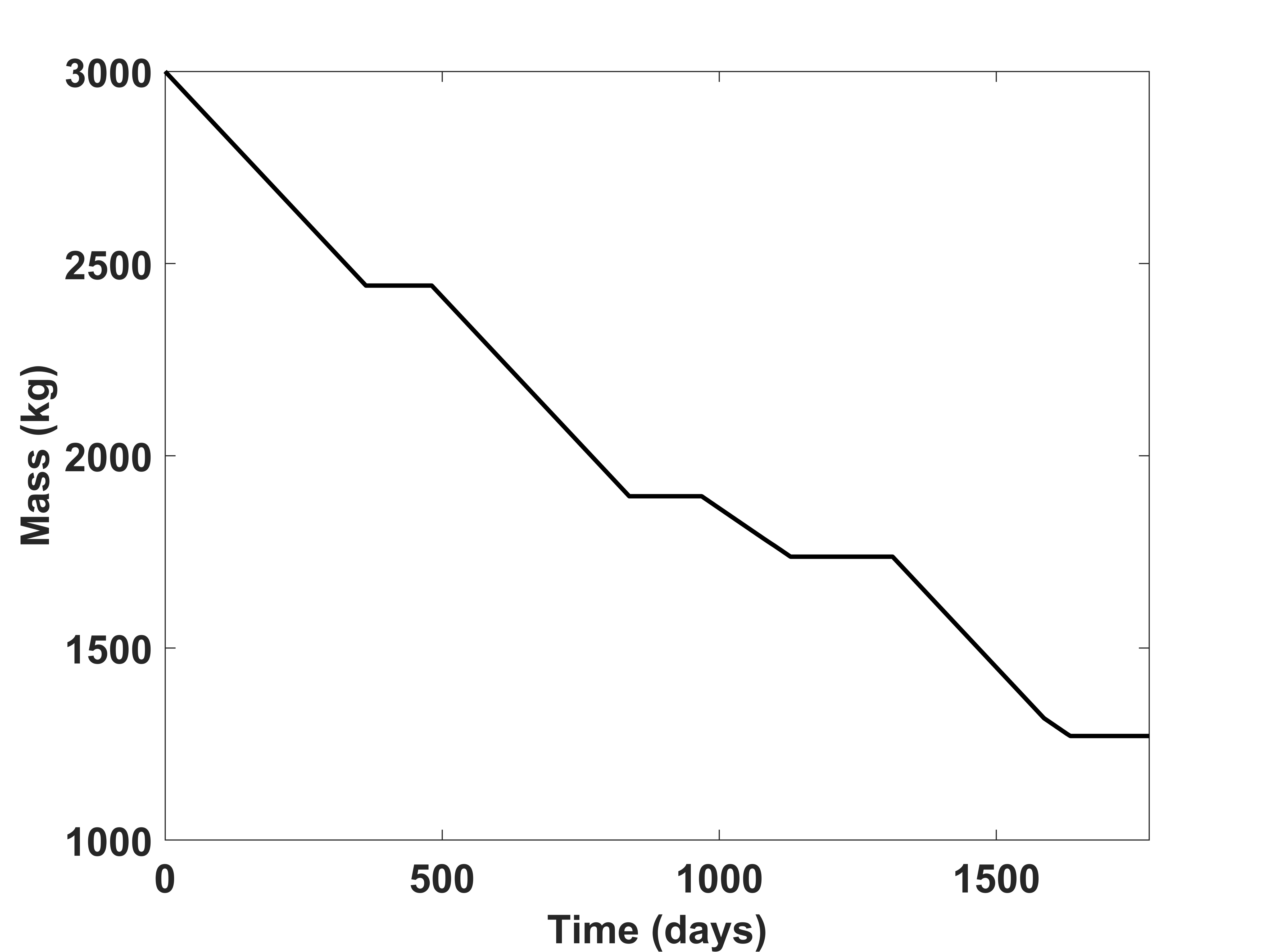}
  \caption{}
  \label{fig:mass_2mode}
\end{subfigure}
\begin{subfigure}{0.5\textwidth}
  \centering
\includegraphics[width=1\linewidth]{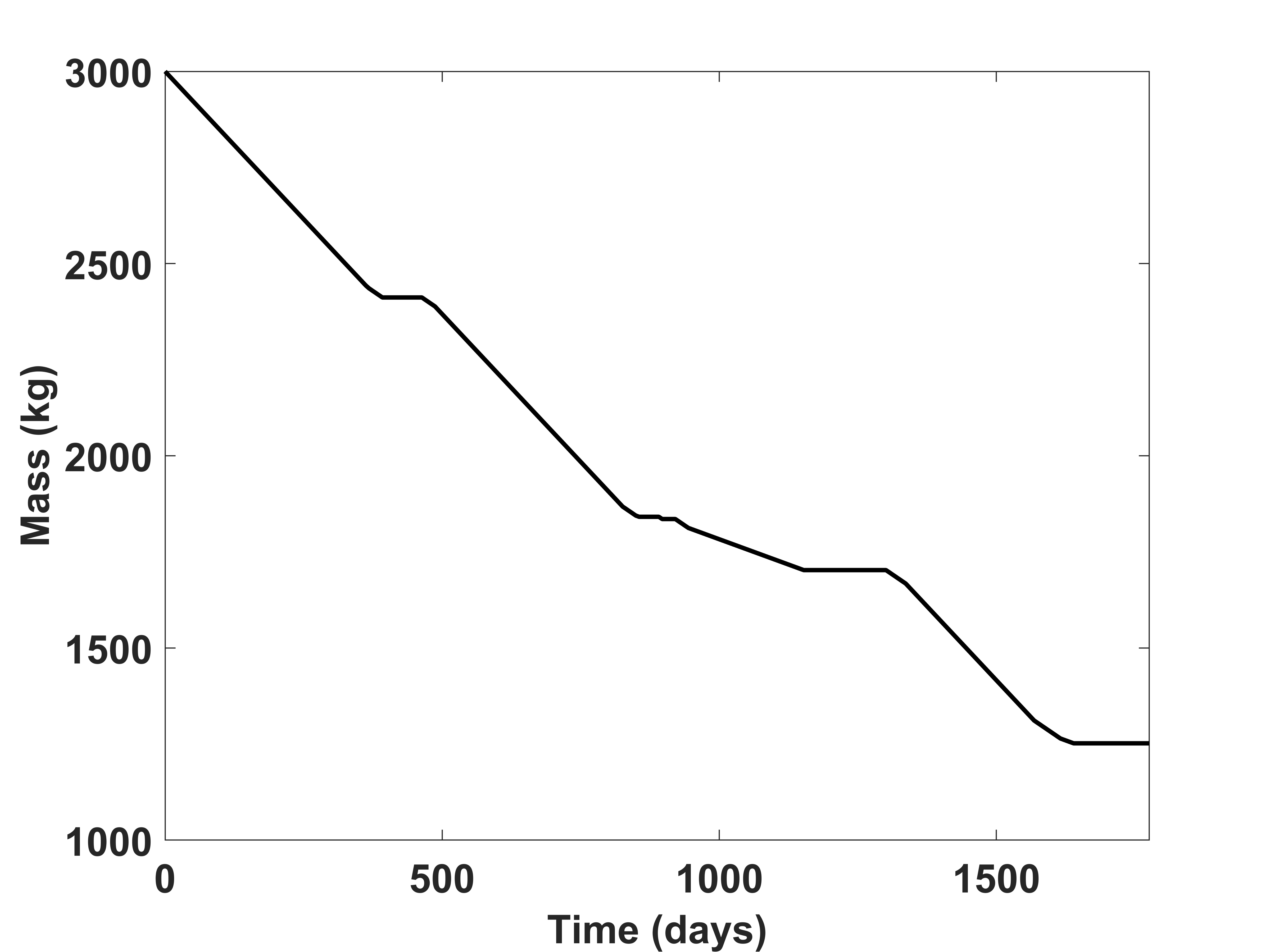}
  \caption{}
  \label{fig:mass_3mode}
\end{subfigure}%
\begin{subfigure}{0.5\textwidth}
  \centering
\includegraphics[width=1\linewidth]{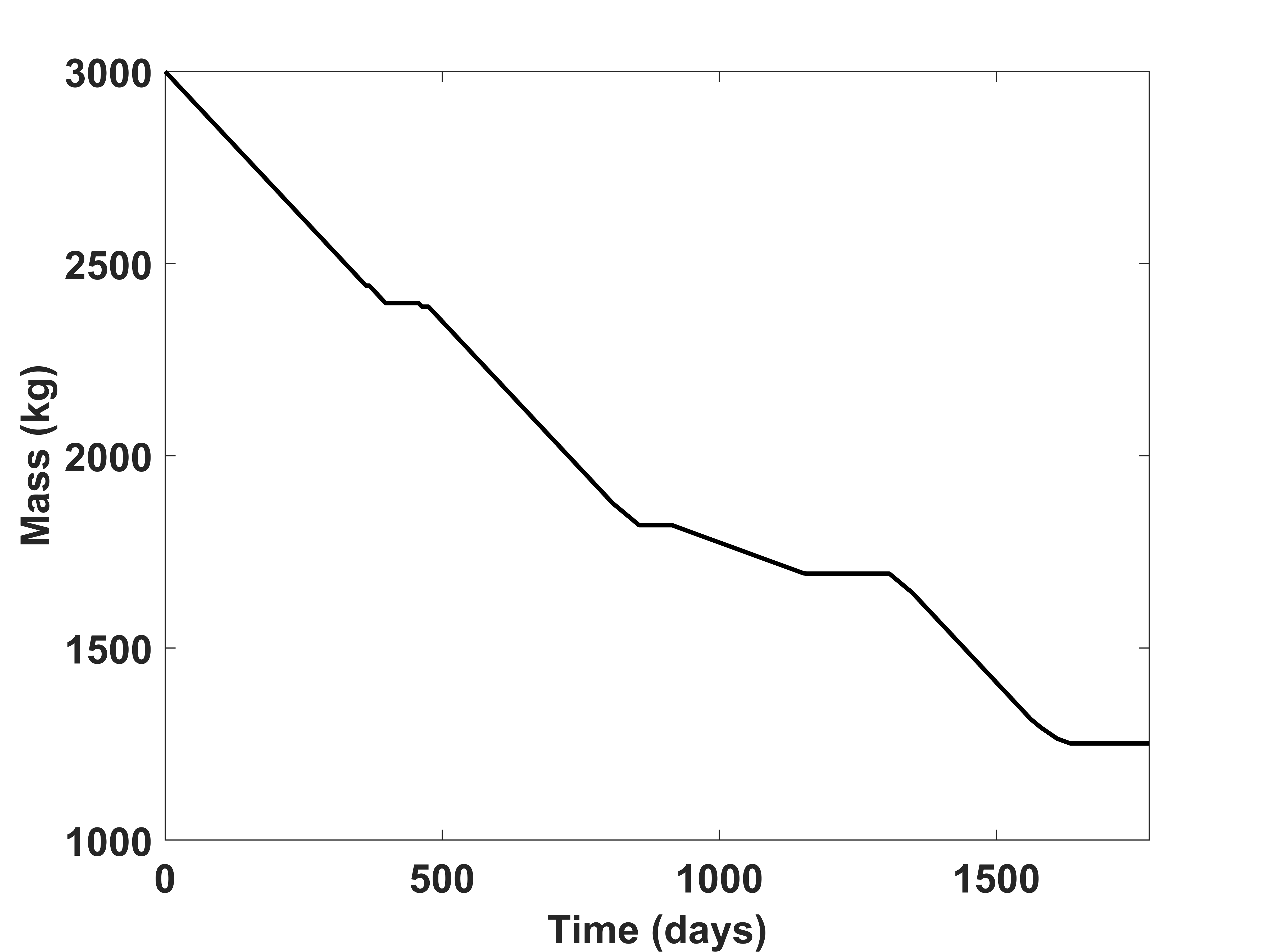}
  \caption{}
  \label{fig:mass_4mode}
\end{subfigure}

\caption{Earth-67P problem total mass changes vs. time: (a) one mode, (b) two modes, (c) three modes, and (d) four modes.}
\label{fig:mass}
\end{figure}

 Obtained trajectories are shown in Figure~\ref{fig:traj}, the thrust profiles are shown in Figure~\ref{fig:thrust}, the power profiles are shown in Figure~\ref{fig:power}, activation vectors are depicted in Figure~\ref{fig:eta} and the fuel consumption profiles are shown in Figure~\ref{fig:mass} for each different number of operation mode case.

\begin{table}[htpb]
    \centering
    \caption{Comparison of optimal solutions: mass breakdown and $P_\text{BL}$ for each subproblem.}
    \begin{tabular}{ccccc}
    \hline
       $N_{mode}$  & 1 & 2 & 3 & 4 \\ \hline
        $m_u$ (kg) & $819.6102$ & $869.5185$ & $868.6574$ & $874.5728$\\
        $m_f$ (kg)& $1238.2003$ & $1270.9570$ & $1251.8866$ & $1251.5264$\\
        $m_\text{SA}$ (kg) & 169.4651& 155.5892& 135.4729 & 129.1613\\
        $m_\text{PSPU}$ (kg) & 242.4101 & 228.5342& 208.4179& 202.1063\\
        $m_\text{PSFS}$ (kg) & 1937.9797 & 1901.9473& 1922.9247& 1923.3209\\
        $P_\text{BL}$ (W) & $16946.5070$ & $15558.9256$ & $13547.2908$ & $12916.1269$\\
        Opt. Modes & $3$ & $[3,20]$  & $[3, 20, 21]$ & $[3, 11, 20, 21]$\\
        $\rho_p$, $\rho_e$& $8.85 \times 10^{-4}$, $8.85 \times 10^{-4}$ & $10^{-4}$, $10^{-4}$& $3 \times 10^{-4}$, $10^{-4}$& $3 \times 10^{-4}$, $10^{-4}$\\
        \hline
    \end{tabular}

    \label{tab:results}
\end{table}

The trajectories are shown in Figure~\ref{fig:traj}. The blue and green dashed orbits represent Earth and Comet 67P, respectively. The coast arcs are shown in blue, while the thrusting arcs are depicted in red. Early arrival/rendezvous (introduced in \cite{taheri2020many}) with Comet 67P is observed in all cases, as a final  coast arc exists at the Comet 67P orbit. There are two main coast arcs in the one-mode case solution, as it is shown in Figure~\ref{fig:traj_1mode}. The spacecraft is on a more inclined orbit until the first coast arc compared to the cases with more modes. Main thrusting arcs are at similar locations for two-, three-, and four-mode cases as they are depicted in Figures~\ref{fig:traj_2mode}, \ref{fig:traj_3mode}, and \ref{fig:traj_4mode}, respectively.

The thrust profiles for each case are shown in Figure~\ref{fig:thrust}, where the discretization nodes are depicted with orange hollow circle markers. For the one-mode case, there is only one thrust level, as shown in Figure~\ref{fig:thrust_1mode}. There is a large thrust arc until approximately 800 days, which increases inclination, as observed in the corresponding trajectory in Figure~\ref{fig:traj_1mode}. 
In Figure~\ref{fig:thrust}, two coast arcs can be observed. For the rest of the cases, four main thrust arcs exist. Around 1000 days, the thrust level switches to a lower thrust mode in Figure~\ref{fig:thrust_2mode} for the two-mode case. The second thrust mode is also observed as the power generated by the solar arrays decreases significantly after 1500 days. In the second mode thrust arc after 1000 days, there is the node that is not exactly at the required thrust level. This can be explained by the fact that $\bm P_\text{E}$ is exactly equal to the $P_\text{sel}$ value of that mode. Therefore, the smooth mode selection introduces a small contribution from the coasting mode. 

In Figure~\ref{fig:thrust_3mode}, three different thrust levels are activated for the three-mode case. At the beginning of the thrust arc where the lowest thrust level (3rd mode) of coasting occurs, there are a couple of nodes that persist in staying at the higher thrust/power level mode (2nd mode). One of the nodes leads to a spike in $\bm P_\text{E}$ during the coast arc just before 100 days, which is also visible in the power profile given in Figure~\ref{fig:power_2mode}. The rest of the nodes correspond to the 2nd mode, since the $\bm P_\text{ava}$ is approximately equal to the $\bm P_\text{E}$ (see Figure~\ref{fig:power_3mode}) and the corresponding power value is slightly higher than the power value of the 2nd mode. Therefore, the 2nd mode is selected as the active mode. If the value of $P_\text{BL}$ is slightly smaller, $\bm P_\text{ava}$ would be slightly smaller so that those nodes could stay at the 3rd mode. This could potentially increase the useful mass of the three-mode case and create a better cost than the two-mode case. 

For the four-mode case, there are four different thrust levels visible in the thrust profile given in Figure~\ref{fig:thrust_4mode}. Unlike the three-mode case, during the thrust arc of around 1000 days, only the 4th mode is active as that is the only available mode in terms of power, which can be seen in Figure~\ref{fig:power_4mode}. We note that the four-mode solution has a lower $P_\text{BL}$ value compared to the three-mode case. As the power gradually decreases after 1500 days (see Figure~\ref{fig:power_4mode}) the thrust level decreases gradually. Around 500 days, a couple of nodes switched from thrust to coast or coast to thrust. By analyzing the power profile given in Figure~\ref{fig:power_4mode}, the $\bm P_\text{E}$ has multiple spikes around that time, which results in power spikes. Also, this coast arc is purely due to the optimality, as the power is already available to activate that mode. These spikes can be attributed to the artifacts of the continuation process performed or they could be due to the flexibility of the optimizer to choose a $\bm P_\text{E}$ that is not constrained by the power availability but only constrained by the optimality of the solution. 

The power profiles of each different case are shown in Figure~\ref{fig:power}. The value of $\bm P_{\text{E}}$ is shown in black color. The purple dashed profile belongs to the selected operating mode powers determined using $\bm P_{\text{E}}$. The blue-colored profile belongs to the power generated by the solar arrays, $\bm P_\text{SA}$. The available power, $\bm P_\text{ava}$, is depicted with red color. One can observe that for each case the $\bm P_\text{E}$ is bounded by $P_\text{max}$ since the upper bound of $\bm P_{\text{E}}$ is $\bm P_{\text{ava}}$ which has an upper bound of $P_\text{max}$. The $\bm P_{\text{ava}}$ has an offset from the $\bm P_\text{SA}$ with the value of $P_\text{sys}$ since $\bm P_{\text{ava}} < P_\text{max}$. Therefore, the power constraints, given in Eq.~\eqref{eq:powconst}, are satisfied. 
Analyzing the one-mode case, given in Figure~\ref{fig:power_1mode}, as the available power decreases, the coasting starts since no other mode is available. Again, just after 1500 days, the available power drops below the mode power's value, resulting in a coast arc. For two-, three- and four-mode cases, the available power decreases around 1000 days leading to a thrust arc with the lowest power mode as shown in Figures~\ref{fig:power_2mode}, \ref{fig:power_3mode}, and \ref{fig:power_4mode}. Similarly, the power decreases after 1500 days, leading to multiple mode switches. The coast arcs around 500 days are due to the optimality since there is enough power available to activate the modes. Coast arcs also happen immediately before and after the thrust arc at 1000 days. In Figure~\ref{fig:power_3mode}, there is a gradual drop in the $\bm P_\text{E}$ around 500 days, unlike the two-mode case. Potentially, directly transitioning to a coast arc would result in a larger coast and could increase the cost. The spike of $\bm P_\text{E}$ around 900 days results in a sudden switch to a thrust arc since there is enough power and this coast is only due to the optimality. In Figure~\ref{fig:power_4mode}, there are spikes in the $\bm P_\text{E}$ value for the four-mode case. This results in sudden changes in the corresponding nodes' value in the thrust profile and the selected power value of the modes. Then, this behavior of the $\bm P_\text{E}$ happens during coast arcs purely from the optimality.

 In Figure~\ref{fig:eta}, the mode switches are shown with $\eta$ values. Each color represents a different mode. It can be seen that at least one mode is activated for each time step. Coasting is considered as another mode with zero thrust value. Therefore, during coast arcs, the coasting mode is active. For the one-mode case, blue is the thrust and orange is the coast, as plotted in Figure~\ref{fig:eta_1mode}. In Figure~\ref{fig:eta_2mode}, three colors represent the modes and the coast. Around 1000 days, the mentioned contribution from the coast and the second mode is visible, with a small drop in its value from 1. Four colors are visible in Figure~\ref{fig:eta_3mode}, corresponding to the three thrust modes and one coast mode. The short spikes in the nodes of $\bm P_\text{E}$ profile correspond to the short-duration activation on $\eta$ profiles. Similarly, in the four-mode case plotted in Figure~\ref{fig:eta_4mode}, five different colors are visible corresponding to four thrust modes and one coast mode. The spikes are around 500 days with short activation times.

The mass consumption profiles for each mode are provided in Figure~\ref{fig:mass}. Whenever a mode change happens, the slope of the mass profile changes as the mass flow rate is changed, as given in Eq.~\eqref{eq:mdot}. The maximum fuel saving corresponds to the two-mode case, as shown in Figure~\ref{fig:mass_2mode} but with a larger solar array compared to the three- and four-mode cases. We note that as the solar array sizing and trajectory optimization are performed simultaneously, the continuation steps are important. As the value of $P_\text{BL}$ potentially changes at each continuation step, a new spacecraft design is optimized for a minimum-fuel expenditure simultaneously. During the continuation process, as the $\rho_p$ decreases, the constraint on the power is enforced incrementally. Therefore, the optimizer can increase the $P_\text{BL}$ to generate more power. The small changes in the $P_\text{BL}$ can lead to a different local optimal solution. Our results indicate that the problem is sensitive to the changes in the value of $P_\text{BL}$. 

\section{Conclusion}
A direct method is proposed to co-optimize spacecraft, trajectory, and solar electric propulsion thruster operation modes. An SPT-140 thruster is considered with discrete operation modes. Each mode is characterized by a specific power input, thrust value, and mass flow rate. The power model includes the capacity of the power processing unit and position- and time-dependent factors capturing changes in power produced by solar arrays. 

The inherent discontinuities due to the power constraint and thruster operation modes are alleviated by leveraging a composite smoothing control (CSC) logic. Sizing of the solar arrays is performed by optimizing the beginning-of-life power of the spacecraft, which appears in the cost function when maximizing the useful mass of the spacecraft is formulated. Upon using the CSC logic, a two-parameter regularized optimization problem is formulated. The regularized optimal control problem is transcribed to a Nonlinear Programming (NLP) problem using the CasADi tool and solved with the IPOPT NLP solver. 

The results show that considering multiple operation modes can increase the useful mass delivered by the spacecraft. In addition, the fuel consumption is also decreased when more modes are considered compared to a single-mode thruster. The beginning-of-life power of the solar arrays decreases as the number of modes increases, which results in lighter solar arrays. Therefore, the co-optimization improves the mass budget of the spacecraft and leads to a better spacecraft design and a more optimal trajectory. The future work includes considering the co-optimization problem with more operation modes and solving different interplanetary missions with different thruster models and comparison of the results against advanced indirect optimization methods. 

\section{Acknowledgment}
The authors would like to acknowledge the NASA Alabama EPSCoR Research Seed Grant program for supporting this research.

\bibliographystyle{AAS_publication}   
\bibliography{references}   
\end{document}